\documentclass[a4paper,11pt]{article}
\usepackage{graphicx}
\usepackage{fancybox}
\usepackage{amssymb}
\usepackage[all]{xy}
\usepackage{color}

\def\rien{\rule{0pt}{0pt}}
\begin{document}
%
%%%%%%%%%%%%%%%%%%%%%%%%%%%%%%%%%%
% Title
%
\title{Learning a Machine for the Decision in a Partially Observable Markov Universe}
%
%%%%%%%%%%%%%%%%%%%%%%%%%%%%%%%%%%
% Author
%
\author{Fr\'ed\'eric Dambreville\\
D\'el\'egation G\'en\'erale pour l'Armement, DGA/CTA/DT/GIP\\
16 Bis, Avenue Prieur de la C\^ote d'Or\\
F 94114, France\\
Email: {\tt Frederic.DAMBREVILLE@dga.defense.gouv.fr}}
%
% make the title area
\maketitle
\begin{abstract}
In this paper, we are interested in optimal decisions in a partially observable Markov universe.
Our viewpoint departs from the dynamic programming viewpoint:
we are directly approximating an optimal strategic tree depending on the observation.
This approximation is made by means of a parameterized probabilistic law.
In this paper, a particular family of hidden Markov models, with input \emph{and} output, is considered as a learning framework.
A method for optimizing the parameters of these HMMs is proposed and applied.
This optimization method is based on the cross-entropic principle.
\end{abstract}
%
% Keywords format
\newcommand{\keywordsname}{Keywords}
\newenvironment{keywords}%
  {\small
    \list{}{\labelwidth0pt
      \leftmargin0pt \rightmargin\leftmargin
      \listparindent\parindent \itemindent0pt
      \parsep0pt
      \let\fullwidthdisplay\relax}%
    \item[\hskip\labelsep\bfseries\keywordsname:]}{\endlist}

\paragraph{Notations.}
Some specific notations are used in this document.
\begin{itemize}
\item The variables $x$, $y$, $z$ and $m$ are used for the action, observation, world state and machine memory,
\item The time $t$ is starting from stage $1$ to the maximal stage $T$.
Variables with subscript outside this scope are synonymous to $\emptyset$.
For example, $\prod_{t=1}^T\pi(x_t|x_{t-1})$ means $\pi(x_1|\emptyset)\prod_{t=2}^T\pi(x_t|x_{t-1})$\,, \emph{ie.} a Markov chain.
A similar principle is used for the \emph{level} supscript $\lambda$ in the definition of hierarchical HMM,
\item Braces used with subscripts, \emph{eg.} $\{\}_{t=1}^T$, only have a grammatical meaning.
More precisely, it means that the symbols inside the braces are duplicated and concatenated according to the subscript.
For example, $\{f_k(\}_{k=1}^3x\{)\}_{k=3}^1$ means $f_1(f_2(f_3(x)))$ and $\{x_k\}_{k=t}^{T}$ means $\prod_{k=t}^{T}x_k$\,,
\item The generic notation for a probability is $P$. However, the functions $p$, $\pi$ and $h$ denote some specific components of the probability.
$p$ is the law of the observation $y$ and state $z$ conditionally to the action $x$.
$\pi$ is a stochastic policy, \emph{ie.} a law of the action conditionally to the observation.
$h$ is an approximation of $\pi$ by a HMM family.
The hidden state of $h$ is defined as the machine memory $m$.
\end{itemize}
\section{Introduction}
There are different degrees of difficulty in planning and control problems.
In most problems, the planner have to start from a given state and terminate in a required final state.
There are several transition rules, which condition the sequence of decision.
For example, a robot may be required to move from room A, starting state, to room B, final state; its action could be \emph{go forward}, \emph{turn right} or \emph{turn left}, and it cannot cross a wall; these are the conditions over the decision.
A first degree in the difficulty is to find at least one solution for the planning.
When the states are only partially known or the actions are not deterministic, the difficulty is quite enhanced:
the planner has to take into account the various observations.
Now, the problem becomes much more complex, when this planning is required to be optimal or near-optimal.
For example, find the shortest trajectory which moves the robot from room A to room B.
There are again different degrees in the difficulty, depending on the problem to be deterministic or not, depending on the model of the future observations.
In the particular case of a Markovian problem with the full observation hypothesis, the dynamic programming principle could be efficiently applied (Markov Decision Process theory/MDP).
This solution has been extended to the case of partial observation (Partially Observable Markov Decision Process/POMDP), but this solution is generally not practicable, owing to the huge dimension of the variables.
\\\\
In this paper, we are interested in optimal planning with partial observation.
A Markovian hypothesis is made.
Our viewpoint departs from the dynamic programming viewpoint:
we are directly approximating an optimal strategic tree depending on the observation.
This approximation is made by means of a parameterized conditional probabilistic law, \emph{ie.} actions depending on observations.
The main problems are the choice of the parameterized laws family and the learning of the optimal parameters.
A particular family of (hierarchical) hidden Markov models, with input \emph{and} output, has been used as the learning framework.
A method for optimizing the parameters of these HHMMs is proposed and applied.
This optimization method is based on the cross-entropic principle\cite{IEEEisda:boer}.
The resulting parameterized law could be seen as a \emph{Virtual Machine} which could generate near-optimal dynamic strategies for any drawing of the problem.
\\\\
The next section introduces some formalism and gives a quick description of the MDP/POMDP problems.
It is recalled that the Dynamic Programming method could solve these problems, but that this solution is intractable for an actual POMDP problem.
It is then proposed a new near-optimal planning method, based on the direct approximation of the optimal decision tree.
The third section introduces the family of Hierarchical Hidden Markov Models.
A particular sub-family of HHMMs is proposed as a candidate for the approximation of decision trees.
The fourth section describes the method for optimizing the parameters of the HHMM, in order to approximate the optimal decision tree for the POMDP problem.
The cross-entropy method is described and applied.
The fifth section gives an example of application.
The paper is then concluded.
\section{MDP and POMDP}
\begin{figure}
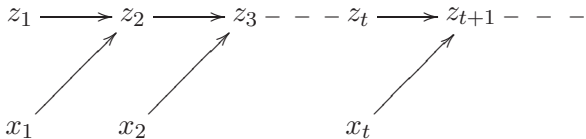

\caption{A Markovian World}
\label{ISDA:fig:1}
\begin{center}
\begin{tabular}{c}
\vspace{-15pt}\\
\xy
(-.5,.5)*+{},
%\green
(0,0)*+{z_1}="b00",
(15,0)*+{z_2}="b01",
(30,0)*+{z_3}="b02",
(45,0)*+{z_t}="b03",
(60,0)*+{z_{t+1}}="b04",
(75,0)="b05",
\ar @{->}"b00";"b01",
\ar @{->}"b01";"b02",
\ar @{--}"b02";"b03",
\ar @{->}"b03";"b04",
\ar @{--}"b04";"b05",\POS "b05",
(0,-15)*+{x_1}="b10",
(15,-15)*+{x_2}="b11",
(45,-15)*+{x_t}="b12",
\ar @{->}"b10";"b01",
\ar @{->}"b11";"b02",
\ar @{->}"b12";"b04",\POS "b04",
(75.5,-15.5)*+{}
\endxy
\vspace{-10pt}
\end{tabular}
\end{center}
\end{figure}

\begin{figure}
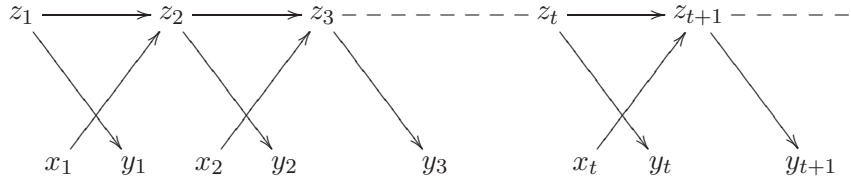

\caption{A Hidden Markov Model (with control)}
\label{ISDA:fig:2}
\begin{center}
\begin{tabular}{c}
\vspace{-15pt}\\
\scalebox{1}{\xy
(-.5,.5)*+{},
%\green
(0,0)*+{z_1}="b00",
(20,0)*+{z_2}="b01",
(40,0)*+{z_3}="b02",
(70,0)*+{z_t}="b03",
(90,0)*+{z_{t+1}}="b04",
(110,0)="b05",
\ar @{->}"b00";"b01",
\ar @{->}"b01";"b02",
\ar @{--}"b02";"b03",
\ar @{->}"b03";"b04",
\ar @{--}"b04";"b05",\POS "b05",
(15,-20)*+{y_1}="b10y",
(35,-20)*+{y_2}="b11y",
(55,-20)*+{y_3}="b12y",
(85,-20)*+{y_t}="b13y",
(105,-20)*+{y_{t+1}}="b14y",
(5,-20)*+{x_1}="b10x",
(25,-20)*+{x_2}="b11x",
(75,-20)*+{x_t}="b12x",
\ar @{->}"b00";"b10y",
\ar @{->}"b01";"b11y",
\ar @{->}"b02";"b12y",
\ar @{->}"b03";"b13y",
\ar @{->}"b04";"b14y",
\ar @{->}"b10x";"b01",
\ar @{->}"b11x";"b02",
\ar @{->}"b12x";"b04",\POS "b04",
(110.5,-20.5)*+{}
\endxy}\vspace{-10pt}
\end{tabular}
\end{center}
\end{figure}

It is assumed that a subject is acting in a given world with a given purpose or mission.
The goal is to optimize the accomplishment of this mission.
As a framework for the optimization, a prior model is hypothetized for the world (world state and observation depending on the action of the subject) and for the mission evaluation (a function depending on world state, observation and action).
In the next paragraphs, a Markovian modelling of the world and a recursive modelling of the evaluation are considered.
\paragraph{The world.}
Our modelling of the world is based on the MDP/POMDP formalism, \emph{ie.} (Partially Observable) Markov Decision Process.
In the MDP or POMDP formalism, the world is considered as a Markov chain conditioned by the past action of the subject.
Assume that the world is characterized by a \emph{state} variable $z$ and that the action is described by a variable $x$.
Assume that the time is starting from $1$, and denote respectively $z_t$ and $x_t$ the world state and the action for the time $t$\,.
The law of $z$ with respect to the past actions $x$ thus verifies the property of Markov:
$$
P(z_t|z_{1:t-1},x_{1:t-1})=p(z_t|z_{t-1},x_{t-1})\;,
$$
where the notation $z_{1:t}$ means $z_1,\dots,z_t$\,, and the notation $z_{0}$\,, \emph{ie.} before starting time, means $\emptyset$\,.
This Markov chain is thus characterized by two functions: an initial law $p(z_1|\emptyset)$ and a transition law $p(z_t|z_{t-1},x_{t-1})$ for $t\ge2$\,.
This distinction is often omitted in this paper, but it is implied.
The law of $z|x$ is represented graphically by the Bayesian Network of figure~\ref{ISDA:fig:1}\,.
In this description, the arrows indicate the dependency between variables.
For example, $z_t\rightarrow z_{t+1}\leftarrow x_t$ means that the law of $z_{t+1}$ is defined conditionally to the past variables $x_t$ and $z_t$\,.
Such Bayesian Network representations will be used several times in this paper for their pedagogic skill.
\paragraph{The observation.}
MDP and POMDP assume that the world is observed during the process.
\vspace{5pt}\\
\emph{In MDP,} the world is fully observed; \emph{the observation at time $t$ is the world state $z_t$\,.}
\vspace{5pt}\\
\emph{In POMDP,} the world is partially observed; \emph{the observation at time $t$ is denoted $y_t$\,.}
The variable of observation $y_t$ is only depending on the current world state $z_t$\,:
$$
P(y_t|z_{1:t},y_{1:t-1},x_{1:t})=p(y_t|z_t)\;.
$$
The variables $y,z|x$ constitute a \emph{Hidden Markov Model} (with control).
This HMM is represented in figure~\ref{ISDA:fig:2}\,.
\paragraph{Evaluation and optimal planning.}
The previous paragraphs have built a modelization of the world, of the actions and of the observations.
We are now giving a characterization of the mission to be accomplished.
\emph{The mission is limited in time.}
Let $T$ be this maximum time.
In the most generality, the mission is evaluated by a function $V(x_{1:T},y_{1:T},z_{1:T})$ defined on the trajectories of $x,y,z$\, ($V(x_{1:T},z_{1:T})$ in the case of a MDP).
Typically, the function $V$ could be used for computing the time needed for the mission accomplishment.
The purpose is to construct an optimal decision tree $x(obs)$ depending on the observation $obs$ in order to maximize the \emph{mean evaluation}.
This is a dynamic optimization problem, since the actions depend on the previous observations.
This problem is expressed differently for MDP and POMDP:
\subparagraph{MDP.} $obs_t=z_t$\,. Optimize $x_t(z_{1:t-1})|_{1\le t\le T}$\,:
$$
\mathrm{Find\quad} x_O\in\arg\max_{x}\sum_{z_{1:T}}V(x_t(z_{1:t-1})|_{1\le t\le T},z_{1:T})
\prod_{t=1}^Tp(z_t|z_{t-1},x_{t-1}(z_{1:t-2}))\;.
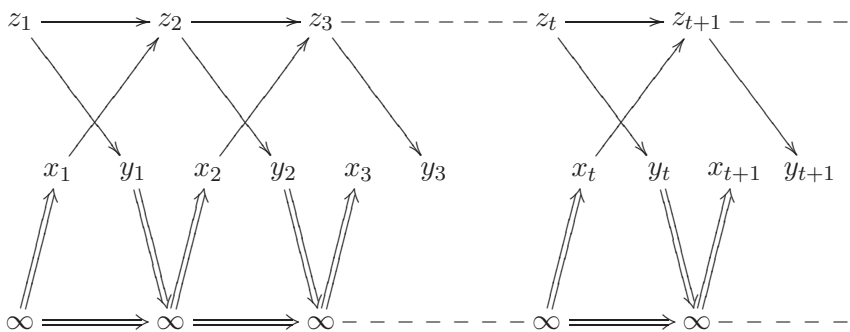
$$
\subparagraph{POMDP.} \rien$obs_t=y_t$\,. Optimize $x_t(y_{1:t-1})|_{1\le t\le T}$\,:
\begin{equation}
\label{Isda:Eq:1}
\begin{array}{@{}l@{}}\displaystyle
\mathrm{Find\quad} x_O\in\arg\max_{x}\sum_{y_{1:T}}\sum_{z_{1:T}}V(x_t(y_{1:t-1})|_{1\le t\le T},y_{1:T},z_{1:T})
\\\displaystyle
\hspace{130pt}\prod_{t=1}^Tp(y_t|z_t)\,p(z_t|z_{t-1},x_{t-1}(y_{1:t-2}))\;.
\end{array}
\end{equation}
\begin{figure}
\caption{MDP planning}
\label{ISDA:fig:3}
\begin{center}
\begin{tabular}{c}
\vspace{-15pt}\\
\xy
(-.5,.5)*+{},
%\green
(0,0)*+{z_1}="b00",
(15,0)*+{z_2}="b01",
(30,0)*+{z_3}="b02",
(45,0)*+{z_t}="b03",
(60,0)*+{z_{t+1}}="b04",
(75,0)="b05",
\ar @{->}"b00";"b01",
\ar @{->}"b01";"b02",
\ar @{--}"b02";"b03",
\ar @{->}"b03";"b04",
\ar @{--}"b04";"b05",\POS "b05",
(0,-15)*+{x_1}="b10",
(15,-15)*+{x_2}="b11",
(30,-15)*+{x_3}="b12",
(45,-15)*+{x_t}="b13",
(60,-15)*+{x_{t+1}}="b14",
\ar @{->}"b10";"b01",
\ar @{->}"b11";"b02",
\ar @{->}"b13";"b04",\POS "b04",
(0,-30)*+{\infty}="b00I",
(15,-30)*+{\infty}="b01I",
(30,-30)*+{\infty}="b02I",
(45,-30)*+{\infty}="b03I",
(60,-30)*+{\infty}="b04I",
(75,-30)="b05I",
\ar @2{-->}"b00I";"b01I",
\ar @2{-->}"b01I";"b02I",
\ar @{--}"b02I";"b03I",
\ar @2{-->}"b03I";"b04I",
\ar @2{->}"b00";"b01I",
\ar @2{->}"b01";"b02I",
\ar @2{->}"b03";"b04I",
\ar @2{->}"b00I";"b10",
\ar @2{->}"b01I";"b11",
\ar @2{->}"b02I";"b12",
\ar @2{->}"b03I";"b13",
\ar @2{->}"b04I";"b14",
\ar @{--}"b04I";"b05I",
\POS "b05I",
(75.5,-30.5)*+{}
\endxy\vspace{-10pt}
\end{tabular}
\end{center}
\end{figure}
\begin{figure}
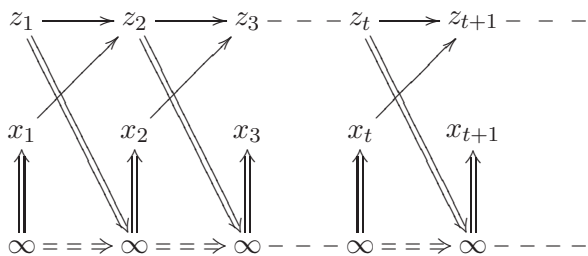

\caption{POMDP planning}
\label{ISDA:fig:4}
\begin{center}
\begin{tabular}{c}
\vspace{-15pt}\\
\scalebox{1}{\xy
(-.5,.5)*+{},
%\green
(0,0)*+{z_1}="b00",
(20,0)*+{z_2}="b01",
(40,0)*+{z_3}="b02",
(70,0)*+{z_t}="b03",
(90,0)*+{z_{t+1}}="b04",
(110,0)="b05",
\ar @{->}"b00";"b01",
\ar @{->}"b01";"b02",
\ar @{--}"b02";"b03",
\ar @{->}"b03";"b04",
\ar @{--}"b04";"b05",\POS "b05",
(15,-20)*+{y_1}="b10y",
(35,-20)*+{y_2}="b11y",
(55,-20)*+{y_3}="b12y",
(85,-20)*+{y_t}="b13y",
(105,-20)*+{y_{t+1}}="b14y",
(5,-20)*+{x_1}="b10x",
(25,-20)*+{x_2}="b11x",
(45,-20)*+{x_3}="b12x",
(75,-20)*+{x_t}="b13x",
(95,-20)*+{x_{t+1}}="b14x",
(0,-40)*+{\infty}="b00I",
(20,-40)*+{\infty}="b01I",
(40,-40)*+{\infty}="b02I",
(70,-40)*+{\infty}="b03I",
(90,-40)*+{\infty}="b04I",
(110,-40)="b05I",
\ar @{->}"b00";"b10y",
\ar @{->}"b01";"b11y",
\ar @{->}"b02";"b12y",
\ar @{->}"b03";"b13y",
\ar @{->}"b04";"b14y",
\ar @{->}"b10x";"b01",
\ar @{->}"b11x";"b02",
\ar @{->}"b13x";"b04",
\POS "b04",
\ar @2{->}"b00I";"b01I",
\ar @2{->}"b01I";"b02I",
\ar @{--}"b02I";"b03I",
\ar @2{->}"b03I";"b04I",
\ar @2{->}"b10y";"b01I",
\ar @2{->}"b11y";"b02I",
\ar @2{->}"b13y";"b04I",
\ar @2{->}"b00I";"b10x",
\ar @2{->}"b01I";"b11x",
\ar @2{->}"b02I";"b12x",
\ar @2{->}"b03I";"b13x",
\ar @2{->}"b04I";"b14x",
\ar @{--}"b04I";"b05I",
\POS "b05I",
(110.5,-40.5)*+{}
\endxy}\vspace{-10pt}
\end{tabular}
\end{center}
\end{figure}
\\[0pt]These optimizations are schematized in figures~\ref{ISDA:fig:3} and~\ref{ISDA:fig:4}.
In these figures, the double arrows $\Rightarrow$ characterize the variables to be optimized.
More precisely, these arrows describe the flow of information between the observations and the actions.
\emph{The cells denoted $\infty$ are making decisions and transmitting all the received and generated information} (including the actions).
This architecture illustrates that planning with observation is an \emph{in(de)finite-memory} problem\,: the decision depends on the whole past observations.
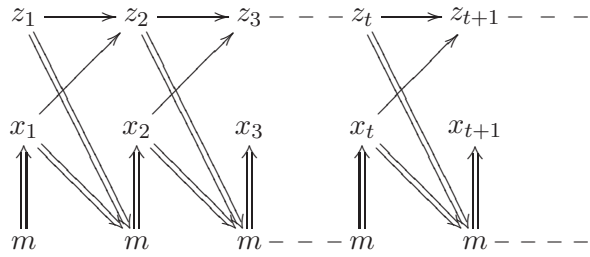
But in the particular case of MDP, the problem is finite memory: the dashed arrows of figure~\ref{ISDA:fig:3} may be removed resulting in the finite memory BN of figure~\ref{ISDA:fig:3b} (memory of the past decision is maintained).
\vspace{5pt}\\
This ``finiteness'' is an ingredient which makes the MDP problems rather easily solvable by means of the dynamic programming method.
DP methods are introduced in section~\ref{Isda:DPmethod}.
As a result, DP produces the root $x_{O,1}$ of the optimal decision tree, which is a sufficient tool for the decision in an actual open-loop planning process.
More precisely, the whole process for optimizing the actual mission takes the form:
\begin{enumerate}
\item Initialize the MDP/POMDP priors of the mission, 
\item \label{Isda:Algo:1:step:1} Compute the theoretical optimal decision root $x_{O,1}$,
\item The decision $x_{O,1}$ is applied to the actual world,
\item Make actual observations and update the prior; in particular, the time is forwarded by $1$,
\item Repeat from step~\ref{Isda:Algo:1:step:1} until accomplishment of the mission.
\end{enumerate}
\begin{figure}
\caption{MDP planning with finite memory}
\label{ISDA:fig:3b}
\begin{center}
\begin{tabular}{c}
\vspace{-15pt}\\
\xy
(-.5,.5)*+{},
%\green
(0,0)*+{z_1}="b00",
(15,0)*+{z_2}="b01",
(30,0)*+{z_3}="b02",
(45,0)*+{z_t}="b03",
(60,0)*+{z_{t+1}}="b04",
(75,0)="b05",
\ar @{->}"b00";"b01",
\ar @{->}"b01";"b02",
\ar @{--}"b02";"b03",
\ar @{->}"b03";"b04",
\ar @{--}"b04";"b05",\POS "b05",
(0,-15)*+{x_1}="b10",
(15,-15)*+{x_2}="b11",
(30,-15)*+{x_3}="b12",
(45,-15)*+{x_t}="b13",
(60,-15)*+{x_{t+1}}="b14",
\ar @{->}"b10";"b01",
\ar @{->}"b11";"b02",
\ar @{->}"b13";"b04",\POS "b04",
(0,-30)*+{m}="b00I",
(15,-30)*+{m}="b01I",
(30,-30)*+{m}="b02I",
(45,-30)*+{m}="b03I",
(60,-30)*+{m}="b04I",
(75,-30)="b05I",
\ar @2{->}"b10";"b01I",
\ar @2{->}"b11";"b02I",
\ar @2{->}"b13";"b04I",
\ar @{--}"b02I";"b03I",
\ar @2{->}"b00";"b01I",
\ar @2{->}"b01";"b02I",
\ar @2{->}"b03";"b04I",
\ar @2{->}"b00I";"b10",
\ar @2{->}"b01I";"b11",
\ar @2{->}"b02I";"b12",
\ar @2{->}"b03I";"b13",
\ar @2{->}"b04I";"b14",
\ar @{--}"b04I";"b05I",
\POS "b05I",
(75.5,-30.5)*+{}
\endxy\vspace{-10pt}
\end{tabular}
\end{center}
\end{figure}
\subsection{Dynamic Programming method}
\label{Isda:DPmethod}
This presentation is extremely simplified.
More detailed references are available\cite{IEEEisda:bellman}\cite{IEEEisda:sondik}\cite{IEEEisda:rocco}\,.
The Dynamic Programming method works for an additive evaluation $V(x_{1:T},z_{1:T})=\sum_{t=1}^TV_t(x_{t},z_{t})$\,.
In this simplified presentation, a final evaluation is assumed, \emph{ie.} $V(x_{1:T},z_{1:T})=V_T(x_T,z_T)$\,, but the principle is the same.
\paragraph{The MDP case.}
To be computed:
$$
x_{O,1}\in\arg_{x_1}\max_{x}\sum_{z_{1:T}}V_T(x_T(z_{1:T-1}),z_{T})
\prod_{t=1}^Tp(z_t|z_{t-1},x_{t-1}(z_{1:t-2}))\;.
$$
This computation is easily refined:
\begin{equation}
\label{DPM:1}
x_{O,1}\in\arg_{x_1}\left\{\max_{x_t}\sum_{z_t}p(z_t|z_{t-1},x_{t-1})\right\}_{t=1}^TV_T(x_T,z_{T})\;.
\end{equation}
This factorization is deeply related to the finite memory property of the MDP.
More precisely, the ``optimal'' evaluation for stage $t$, denoted $W_t$\,, may be computed recursively, as a direct result of (\ref{DPM:1})\,:
$$
\left\{\begin{array}{@{}l@{}}\displaystyle\vspace{5pt}
W_T(x_T,z_T)=V_T(x_T,z_T)\;,
\\\displaystyle\vspace{5pt}
W_t(x_t,z_t)=\max_{x_{t+1}}\sum_{z_{t+1}}p(z_{t+1}|z_t,x_t)W_{t+1}(x_{t+1},z_{t+1})\;.\vspace{-15pt}
\end{array}\right.\vspace{15pt}
$$
Associated to these evaluations are defined the optimal strategies $x_{O,t}$ for the stage $t$\,:
$$
x_{O,t}\in\arg\max_{x_t}\sum_{z_t}p(z_t|z_{t-1},x_{t-1})W_t(x_t,z_t)\;.
$$
In particular: $x_{O,1}\in\arg\max_{x_1}\sum_{z_1}p(z_1|\emptyset)W_1(x_1,z_1)\;.$\vspace{5pt}\\
It happens that the optimal strategy for stage $t$ is only depending on the previous action and state $x_{t-1},z_{t-1}$\,.
This is exactly what have been intuited from figure~\ref{ISDA:fig:3b} as a finite memory property.
\\\\
This method, known as the dynamic programming, thus relies on the computation of backwards functions defined over the variables $x_t,z_t$\,.
The DP methods have been extended to POMDP problems, but this time the recursive functions are much more intricated and almost intractable.
\paragraph{The POMDP case.}
To be computed:
$$
\begin{array}{@{}l@{}}\displaystyle
x_{O,1}\in\arg_{x_1}\max_{x}\sum_{y_{1:T}}\sum_{z_{1:T}}V_T(x_T(z_{1:T-1}),y_{T},z_{T})
\\\displaystyle
\hspace{150pt}\prod_{t=1}^Tp(y_t|z_t)\,p(z_t|z_{t-1},x_{t-1}(y_{1:t-2}))\;.
\end{array}
$$
The previous factorization works partially:
$$
x_{O,1}\in\arg_{x_1}\left\{\max_{x_t}\sum_{y_t}\right\}_{t=1}^T\;\sum_{z_{1:T}}V_T(x_T,y_{T},z_{T})
\prod_{t=1}^Tp(y_t|z_t)\,p(z_t|z_{t-1},x_{t-1})\;.
$$
But this factorization is incomplete and the previous recursion is not possible anymore.
Still, there is an answer to this problem by means of a dynamic programming method.
But the recursion cannot rely on the variables $x_t,y_t$ (last action and observation) only: these variables contain insufficient informations to predict the future.
It is necessary to transmit the probabilistic \emph{belief} over $z_t$ estimated from the whole past actions and observations.
Denote $b_t(z_t)$ the belief over $z_t$\,.
The solution of the POMDP is constructed recursively:
$$
\left\{\begin{array}{@{}l@{}}\displaystyle\vspace{5pt}
W_T(b_T)=\max_{x_T}\sum_{y_T}\sum_{z_T}b_T(z_T)p(y_T|z_T)V_T(x_T,y_T,z_T)\;;
\\\displaystyle\vspace{5pt}
\beta_{t+1}(z_{t+1}|x_t,y_t,b_t)=\sum_{z_{t}}b_t(z_t)p(y_t|z_t)p(z_{t+1}|z_t,x_t)\;,
\\\displaystyle\vspace{5pt}
W_t(b_t)=\max_{x_t}\sum_{y_t}W_{t+1}(\beta_{t+1}(\cdot|x_t,y_t,b_t))\;;
\\\displaystyle\vspace{5pt}
%\mbox{where }\mbox{where }
\beta_2(z_2|x_1,y_1)=\sum_{z_1}p(z_1|\emptyset)p(y_1|z_1)p(z_{2}|x_1,z_1)\;,
\\\displaystyle
x_{O,1}\in\arg\max_{x_1}\sum_{y_1}W_{2}(\beta_2(\cdot|x_1,y_1))\;.
\end{array}\right.
$$
This recursive construction of $W_t$ is very difficult, since the variable $b_t$ is continuous and has a very high dimension.
For example, if the state $z_t$ is the position of a target in a discrete space of $20\times 20$ cells (the same space as in the example\footnote{In this example the dimension is even more, \emph{ie.} $4^2400^3-1$\,.} of section~\ref{Isda:example})\,, this dimension is $399$\,.
Of course, there are many possible refinements and approximations, but even so, the DP method becomes easily intractable.
\vspace{5pt}\\
This particular difficulty for solving POMDP problems was illustrated by the figures~\ref{ISDA:fig:3} and~\ref{ISDA:fig:4}\,: MDPs are fundamentally more simple than POMDPs.
Indeed, since the world state $z$ is a Markov chain, the knowledge of $z_t$ together with the last action $x_t$ is sufficient to predict all the future (figure~\ref{ISDA:fig:3b}).
On the other hand, the knowledge of $y_t$ and of $x_t$ is not sufficient to predict the future of the Markov chain $z$.
\subsection{Direct approximation of the decision tree}
\begin{figure}
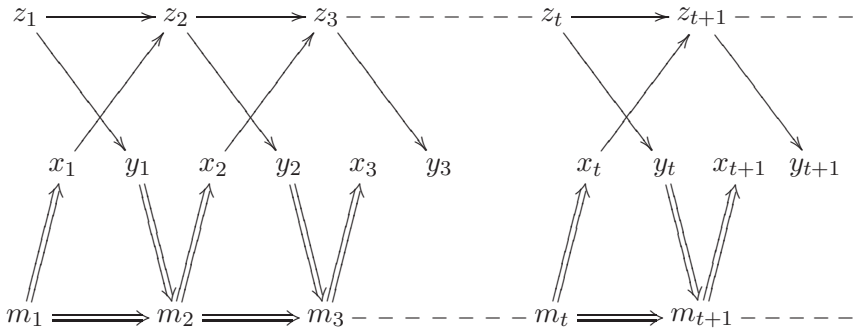

\caption{Finite-memory planning approximation}
\label{ISDA:fig:5}
\begin{center}
\begin{tabular}{c}
\vspace{-15pt}\\
\scalebox{1}{\xy
(-.5,.5)*+{},
%\green
(0,0)*+{z_1}="b00",
(20,0)*+{z_2}="b01",
(40,0)*+{z_3}="b02",
(70,0)*+{z_t}="b03",
(90,0)*+{z_{t+1}}="b04",
(110,0)="b05",
\ar @{->}"b00";"b01",
\ar @{->}"b01";"b02",
\ar @{--}"b02";"b03",
\ar @{->}"b03";"b04",
\ar @{--}"b04";"b05",\POS "b05",
(15,-20)*+{y_1}="b10y",
(35,-20)*+{y_2}="b11y",
(55,-20)*+{y_3}="b12y",
(85,-20)*+{y_t}="b13y",
(105,-20)*+{y_{t+1}}="b14y",
(5,-20)*+{x_1}="b10x",
(25,-20)*+{x_2}="b11x",
(45,-20)*+{x_3}="b12x",
(75,-20)*+{x_t}="b13x",
(95,-20)*+{x_{t+1}}="b14x",
(0,-40)*+{m_1}="b00I",
(20,-40)*+{m_2}="b01I",
(40,-40)*+{m_3}="b02I",
(70,-40)*+{m_{t}}="b03I",
(90,-40)*+{m_{t+1}}="b04I",
(110,-40)="b05I",
\ar @{->}"b00";"b10y",
\ar @{->}"b01";"b11y",
\ar @{->}"b02";"b12y",
\ar @{->}"b03";"b13y",
\ar @{->}"b04";"b14y",
\ar @{->}"b10x";"b01",
\ar @{->}"b11x";"b02",
\ar @{->}"b13x";"b04",
\POS "b04",
\ar @2{->}"b00I";"b01I",
\ar @2{->}"b01I";"b02I",
\ar @{--}"b02I";"b03I",
\ar @2{->}"b03I";"b04I",
\ar @2{->}"b10y";"b01I",
\ar @2{->}"b11y";"b02I",
\ar @2{->}"b13y";"b04I",
\ar @2{->}"b00I";"b10x",
\ar @2{->}"b01I";"b11x",
\ar @2{->}"b02I";"b12x",
\ar @2{->}"b03I";"b13x",
\ar @2{->}"b04I";"b14x",
\ar @{--}"b04I";"b05I",
\POS "b05I",
(110.5,-40.5)*+{}
\endxy}\vspace{-10pt}
\end{tabular}
\end{center}
\end{figure}

While the classical DP method generally fails solving the POMDP optimization, some authors have investigated approximations of the solution.
Our work investigates the direct approximation of the decision tree.
\vspace{5pt}\\
In an optimization problem like~(\ref{Isda:Eq:1})\,, the value to be optimized, $x_O$\,, is a \emph{deterministic} object.
In this precise case, $x_O$ is a tree of decision, that is a function which maps to a decision $x_t$ from any sequence of observation $y_{1:t-1}$\,.
It is possible however to have a probabilistic viewpoint.\footnote{Such probabilistic viewpoint is more common, in fact necessary, in game problems.}
Then the problem is equivalent to finding $x,y\mapsto\pi(x|y)$\,, a probabilistic law of actions conditionally to the \emph{past} observations, which maximizes the mean evaluation:
$$
\begin{array}{@{}l@{}}\displaystyle\vspace{5pt}
V(\pi)=\sum_{x_{1:T}}\sum_{y_{1:T}}\sum_{z_{1:T}}\;\;\prod_{t=1}^T\pi(x_t|x_{1:t-1},y_{1:t-1})
\\\displaystyle
\hspace{100pt}\times P(y_{1:T},z_{1:T}|x_{1:T})\;V(x_{1:T},y_{1:T},z_{1:T})\;,
\end{array}
$$
where:
$$
P(y_{1:T},z_{1:T}|x_{1:T})=\prod_{t=1}^Tp(y_t|z_t)\,p(z_t|z_{t-1},x_{t-1})\;.
$$
If the solution is optimal, there will not be a great difference with the deterministic case:
when the solution $x_O$ is unique, the optimal law $\pi_O$ is a dirac on $x_O$\,.
But things change when approximating $\pi$\,.
Now, why using a probability to approximate the optimal strategy?
The main point is that probabilistic models seem more suitable for approximation, than for example deterministic decision trees.
The second point is that we are sure to approximate continuously: indeed, $\pi\mapsto V(\pi)$ is continuous.
There is a third point, which is more ``philosophic''.
Considering the figure~\ref{ISDA:fig:4}\,, it appears clearly that it is symmetric: the arrows $\Rightarrow$ and $\rightarrow$ play a quite similar role.
In the deterministic viewpoint, this is just a coincidence.
But if we are considering probabilistic strategies, the arrows $\Rightarrow$ in the Bayesian Network of figure~\ref{ISDA:fig:4} are constituting a HMM controlled by the observation, and with infinite memory.
Of course, a natural approximation of a HMM with infinite memory is a HMM with finite memory!
Then replacing the infinite-memory states $\infty$ by a finite-memory variable $m$\,, a perfectly symmetric problem is obtained in figure~\ref{ISDA:fig:5}\,; perhaps some kind of duality relation.
Well, aestheticism is not the unique interest of such approximation.
The most interesting point is that there are many methods which apply on HMM.
Particularly for learning and optimizing parameters.
Then, the approach developped in this paper is quite general and can be split up into two points:
\begin{itemize}
\item Define a family of parameterized HMMs $\mathcal{H}$\,,
\item Optimize the parameters of the HMM in order to maximize the mean evaluation:
$$
\mathrm{Find }\quad h_O\in\arg\max_{h\in\mathcal{H}}V(h)\;.
$$
\end{itemize}
The first point is discussed in the next section, where a sub-class of hierarchical HMMs is defined.
The second point is discussed in the following section, which explains a cross-entropic method for optimizing the parameters.
\section{Hierarchical Hidden Markov Model}
In this work, a significant investment was dedicated to the construction of a ``good'' parameterized probabilitic family.
In order to investigate future complex problems, our choice has focused on hierarchical HMM models.
These models are inspired from biology: to solve a complex problem, factorize it and make decisions in a hierarchical fashion.
Low hierarchies manipulate low level informations and actions, making short-term decisions.
High hierarchies manipulate high level informations and actions (uncertainty is less), making long-term decisions.
\vspace{5pt}\\
The implementation of such hierarchical models have not been completely investigated yet;
a rather simple HHMM has been implemented, and our example tests just needed few hierarchies.
Still, this section is deliberately discussing about HHMM with more hindsight.
\subsection{Generality}
\begin{figure}
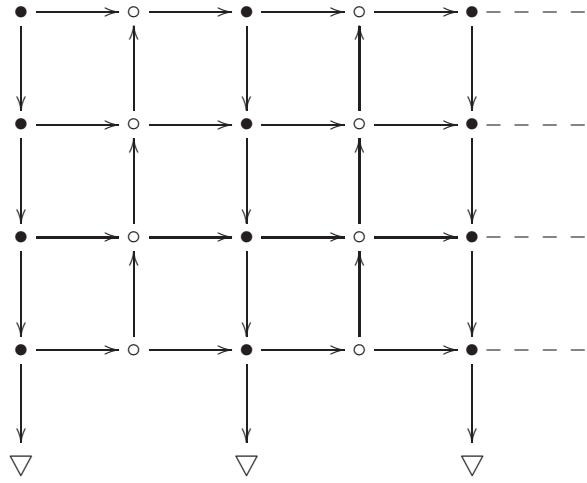

\caption{Model of a Hierarchical HMM}
\label{ISDA:fig:6}
\begin{center}
\begin{tabular}{@{}c@{}}
\vspace{-15pt}\\
\xy
(-.5,.5)*+{},
%\green
(0,0)*+{\bullet}="b00",
(15,0)*+{\circ}="b01",
(30,0)*+{\bullet}="b02",
(45,0)*+{\circ}="b03",
(60,0)*+{\bullet}="b04",
(75,0)="b05",
\ar @{->}"b00";"b01",
\ar @{->}"b01";"b02",
\ar @{->}"b02";"b03",
\ar @{->}"b03";"b04",
\ar @{--}"b04";"b05",\POS "b05",
(0,-15)*+{\bullet}="b10",
(15,-15)*+{\circ}="b11",
(30,-15)*+{\bullet}="b12",
(45,-15)*+{\circ}="b13",
(60,-15)*+{\bullet}="b14",
(75,-15)="b15",
\ar @{->}"b00";"b10",
\ar @{<-}"b01";"b11",
\ar @{->}"b02";"b12",
\ar @{<-}"b03";"b13",
\ar @{->}"b04";"b14",
\ar @{->}"b10";"b11",
\ar @{->}"b11";"b12",
\ar @{->}"b12";"b13",
\ar @{->}"b13";"b14",
\ar @{--}"b14";"b15",\POS "b15",
(0,-30)*+{\bullet}="b00",
(15,-30)*+{\circ}="b01",
(30,-30)*+{\bullet}="b02",
(45,-30)*+{\circ}="b03",
(60,-30)*+{\bullet}="b04",
(75,-30)="b05",
\ar @{<-}"b00";"b10",
\ar @{->}"b01";"b11",
\ar @{<-}"b02";"b12",
\ar @{->}"b03";"b13",
\ar @{<-}"b04";"b14",
\ar @{->}"b00";"b01",
\ar @{->}"b01";"b02",
\ar @{->}"b02";"b03",
\ar @{->}"b03";"b04",
\ar @{--}"b04";"b05",\POS "b05",
(0,-45)*+{\bullet}="b10",
(15,-45)*+{\circ}="b11",
(30,-45)*+{\bullet}="b12",
(45,-45)*+{\circ}="b13",
(60,-45)*+{\bullet}="b14",
(75,-45)="b15",
\ar @{->}"b00";"b10",
\ar @{<-}"b01";"b11",
\ar @{->}"b02";"b12",
\ar @{<-}"b03";"b13",
\ar @{->}"b04";"b14",
\ar @{->}"b10";"b11",
\ar @{->}"b11";"b12",
\ar @{->}"b12";"b13",
\ar @{->}"b13";"b14",
\ar @{--}"b14";"b15",\POS "b15",
(0,-60)*+{\bigtriangledown}="b00",
(30,-60)*+{\bigtriangledown}="b02",
(60,-60)*+{\bigtriangledown}="b04",
\ar @{<-}"b00";"b10",
\ar @{<-}"b02";"b12",
\ar @{<-}"b04";"b14",\POS "b14",
(75.5,-60.5)*+{}
\endxy\vspace{10pt}\\
$\bullet\,/\,\circ\,/\,\bigtriangledown\;=\;$information\,/\,information+boolean\,/\,ouput\vspace{-10pt}
\end{tabular}
\end{center}
\end{figure}

This section gives some general references about HHMMs \emph{without control}.
%The author is not informed of works about controlled HHMM.
However, HHMM have been used yet in control problems.
In \cite{IEEEisda:theo} a HHMM modelling of the world is used in a POMDP problem, but the control (\emph{ie.} the action) is not propagated within the HHMM.
In our models, actions and observations are both propagated within the same HHMM structure.
Thus, the model of controlled HHMM defined here is not related to any existing work, the author knows.
But it is interesting to introduce this section with some results about uncontrolled HHMM.
\\\\
A hierarchical hidden Markov model (HHMM) may be defined as a HMM which output is either a hierarchical HMM or an actual output.
A HHMM could also be considered as a hierarchy of \emph{stochastic} processes calling sub-processes.
It is noteworthy that the length of a sub-process call depends on the sub-process law.
This length is variable.
\vspace{5pt}\\
Simple HMM have been variously applied, for example in speech or hand-writing recognition, in robotic, etc.
However there are still few applications of hierarchical HMM\cite{IEEEisda:theo}\cite{IEEEisda:fine}, although it is probable that HHMM should allow a more abstract representation of the processes.
\\\\
In \cite{IEEEisda:fine}, \emph{Fine, Singer and Tishby} applied HHMM to \emph{identify combinations of letters in handwriting}.
The main difficulty was to derive appropriate Viterbi and Baum-Welch algorithms.
The Complexity of \emph{Fine and al} methods was about $T^3Q^\Lambda$, where $T$ is the length of the samples, $\Lambda$ is the number of \emph{levels} of the HHMM tree and $Q$ is the number of states per levels of the HHMM.
Such required computation times are obviously difficult to apply; for most samples $T$ is too big.
\vspace{5pt}\\
In their work \cite{IEEEisda:murphy}, \emph{Murphy and Paskin} have shown how HHMM could be interpreted as a particular $2-$dimension dynamic Bayesian Network sized $\Lambda\times T$, and derived algorithms of complexity $TQ^{2\Lambda}$.
This is much better, but remains exponential with $\Lambda$.
The number of levels is thus strongly limited.
Although \emph{Murphy and Paskin} did not express this property explicitely, it is easy to show that their Bayesian Network has the Markov property both in time and in level.
In fact, a hierarchical HMM could be interpreted by a Bayesian Network as described in figure~\ref{ISDA:fig:6} (refer to appendix~\ref{sect:App:hhmm})\,.
A hierarchical HMM is thus clearly a HMM with vectorial states.
However, these states should contains some boolean components which are used to define the sub-process call and return\cite{IEEEisda:murphy} (also refer to appendix~\ref{sect:App:hhmm}).
Moreover, there is a downward \emph{and upward} exchange of information between the levels, in a hierarchical HMM.
\vspace{5pt}\\
The next section introduces a model of controlled HHMM\,, which is (freely) inspired from the BN of figure~\ref{ISDA:fig:6}\,.
\begin{figure}
\caption{Model of a controlled Hierarchical HMM}
\label{ISDA:fig:6:bis}
\begin{center}
\begin{tabular}{@{}c@{}}
\vspace{-15pt}\\
\xy
(-.5,.5)*+{},
%\green
(0,0)*+{\circ}="b00",
(15,0)*+{\circ}="b01",
(30,0)*+{\circ}="b02",
(45,0)*+{\circ}="b03",
(60,0)*+{\circ}="b04",
(75,0)="b05",
\ar @{->}"b00";"b01",
\ar @{->}"b01";"b02",
\ar @{->}"b02";"b03",
\ar @{->}"b03";"b04",
\ar @{--}"b04";"b05",\POS "b05",
(0,-15)*+{\circ}="b10",
(15,-15)*+{\circ}="b11",
(30,-15)*+{\circ}="b12",
(45,-15)*+{\circ}="b13",
(60,-15)*+{\circ}="b14",
(75,-15)="b15",
\ar @{->}"b00";"b10",
\ar @{<-}"b01";"b11",
\ar @{->}"b02";"b12",
\ar @{<-}"b03";"b13",
\ar @{->}"b04";"b14",
\ar @{->}"b10";"b11",
\ar @{->}"b11";"b12",
\ar @{->}"b12";"b13",
\ar @{->}"b13";"b14",
\ar @{--}"b14";"b15",\POS "b15",
(0,-30)*+{\circ}="b00",
(15,-30)*+{\circ}="b01",
(30,-30)*+{\circ}="b02",
(45,-30)*+{\circ}="b03",
(60,-30)*+{\circ}="b04",
(75,-30)="b05",
\ar @{<-}"b00";"b10",
\ar @{->}"b01";"b11",
\ar @{<-}"b02";"b12",
\ar @{->}"b03";"b13",
\ar @{<-}"b04";"b14",
\ar @{->}"b00";"b01",
\ar @{->}"b01";"b02",
\ar @{->}"b02";"b03",
\ar @{->}"b03";"b04",
\ar @{--}"b04";"b05",\POS "b05",
(0,-45)*+{\circ}="b10",
(15,-45)*+{\circ}="b11",
(30,-45)*+{\circ}="b12",
(45,-45)*+{\circ}="b13",
(60,-45)*+{\circ}="b14",
(75,-45)="b15",
\ar @{->}"b00";"b10",
\ar @{<-}"b01";"b11",
\ar @{->}"b02";"b12",
\ar @{<-}"b03";"b13",
\ar @{->}"b04";"b14",
\ar @{->}"b10";"b11",
\ar @{->}"b11";"b12",
\ar @{->}"b12";"b13",
\ar @{->}"b13";"b14",
\ar @{--}"b14";"b15",\POS "b15",
(0,-60)*+{\bigtriangledown}="b00",
(15,-60)*+{\bigtriangleup}="b01",
(30,-60)*+{\bigtriangledown}="b02",
(45,-60)*+{\bigtriangleup}="b03",
(60,-60)*+{\bigtriangledown}="b04",
\ar @{<-}"b00";"b10",
\ar @{->}"b01";"b11",
\ar @{<-}"b02";"b12",
\ar @{->}"b03";"b13",
\ar @{<-}"b04";"b14",\POS "b14",
(75.5,-60.5)*+{}
\endxy\vspace{10pt}\\
$\circ\,/\,\bigtriangledown\,/\,\bigtriangleup\;=\;$information+boolean\,/\,ouput\,/\,input\vspace{-10pt}
\end{tabular}
\end{center}
\end{figure}
\subsection{Definition of a controlled model}
\begin{figure}
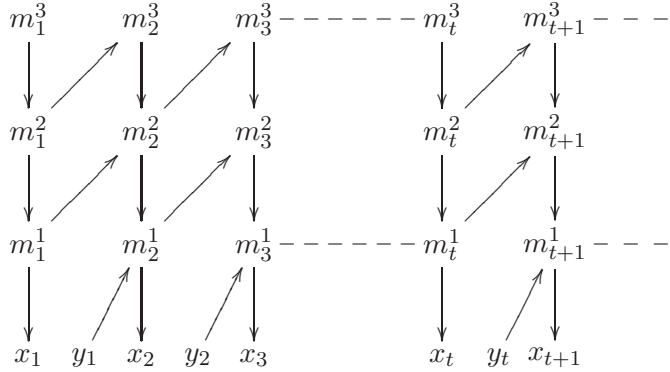

\caption{HHMM model for the planning}
\label{ISDA:fig:7}
\begin{center}
\begin{tabular}{@{}c@{}}
\rien\vspace{-15pt}\\
\xy
(-.5,.5)*+{},
%\green
(0,0)*+{m^3_1}="b00",
(15,0)*+{m^3_2}="b01",
(30,0)*+{m^3_3}="b02",
(55,0)*+{m^3_t}="b03",
(70,0)*+{m^3_{t+1}}="b04",
(85,0)="b05",
(0,-15)*+{m^2_1}="b10",
(15,-15)*+{m^2_2}="b11",
(30,-15)*+{m^2_3}="b12",
(55,-15)*+{m^2_t}="b13",
(70,-15)*+{m^2_{t+1}}="b14",
(85,-15)="b15",
\ar @{->}"b00";"b10",
\ar @{->}"b01";"b11",
\ar @{->}"b02";"b12",
\ar @{->}"b03";"b13",
\ar @{->}"b04";"b14",
\ar @{->}"b10";"b01",
\ar @{->}"b11";"b02",
\ar @{->}"b13";"b04",
\ar @{--}"b02";"b03",
\ar @{--}"b04";"b05",\POS "b05",
(0,-30)*+{m^1_1}="b00",
(15,-30)*+{m^1_2}="b01",
(30,-30)*+{m^1_3}="b02",
(55,-30)*+{m^1_t}="b03",
(70,-30)*+{m^1_{t+1}}="b04",
(85,-30)="b05",
\ar @{<-}"b00";"b10",
\ar @{<-}"b01";"b11",
\ar @{<-}"b02";"b12",
\ar @{<-}"b03";"b13",
\ar @{<-}"b04";"b14",
\ar @{->}"b00";"b11",
\ar @{->}"b01";"b12",
\ar @{->}"b03";"b14",
\ar @{--}"b02";"b03",
\ar @{--}"b04";"b05",\POS "b05",
(0,-45)*+{x_1}="b10",
(15,-45)*+{x_2}="b11",
(30,-45)*+{x_3}="b12",
(55,-45)*+{x_t}="b13",
(70,-45)*+{x_{t+1}}="b14",
(85,-45)="b15",
(7.5,-45)*+{y_1}="b10b",
(22.5,-45)*+{y_2}="b11b",
(62.5,-45)*+{y_t}="b13b",
\ar @{->}"b00";"b10",
\ar @{->}"b01";"b11",
\ar @{->}"b02";"b12",
\ar @{->}"b03";"b13",
\ar @{->}"b04";"b14",
\ar @{->}"b10b";"b01",
\ar @{->}"b11b";"b02",
\ar @{->}"b13b";"b04",\POS "b04",
(85.5,-45.5)*+{}
\endxy\vspace{-10pt}
\end{tabular}
\end{center}
\end{figure}

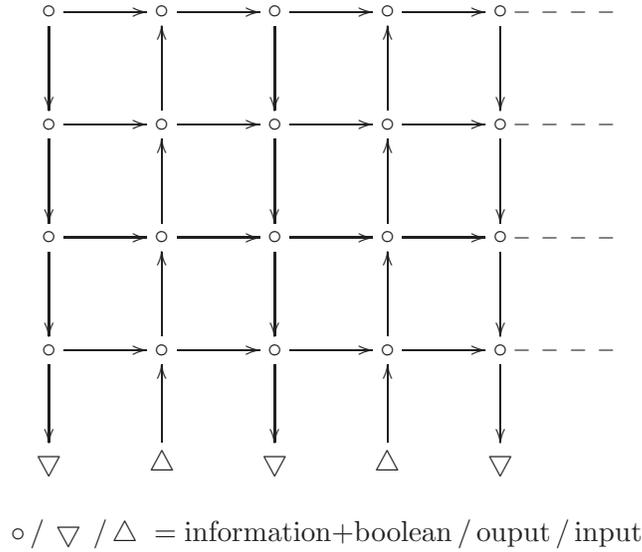
A natural extension of the BN of figure~\ref{ISDA:fig:6} is obtained by adding an input door to each upward column of the BN (then completing it), as described in figure~\ref{ISDA:fig:6:bis}; moreover, boolean components are then needed in all hidden cells, so that the whole memory is required to be discrete or semi-continuous.
Then, an alternance of output and input is obtained, which exactly matches our problem.
This is probably a very good model and future works should investigate this solution.
But for historical reasons in our work, another model, similar, has been chosen for now\dots
\vspace{5pt}\\
In the HHMM family which has been implemented, each memory state receives an information from the current upper-level state and the previous lower-level state.
This model guarantees a propagation of the information between the several levels of the hierarchy, but this propagation is slower than in the case of figure~\ref{ISDA:fig:6:bis}.
More precisely, this HHMM family being denoted $\mathcal{H}$\,, any hhmm $h\in\mathcal{H}$ takes the form:
$$
h(x|y)=\sum_{m_{1:T}^{1:\Lambda}\in M^{1:\Lambda}}\;\prod_{t=1}^Th^0(x_t|m^1_{t})h^1(m_t^1|y_{t-1},m^2_{t})
\prod_{\lambda=2}^{\Lambda}h^\lambda(m_{t}^\lambda|m_{t-1}^{\lambda-1},m_{t}^{\lambda+1})\;,
$$
where $m^\lambda$ is the variable, or memory, for the hidden level $\lambda$ of the HHMM, $M^\lambda$ is the set of possible states for the variable $m^\lambda$\,, and $\Lambda$ is the number of levels of the HHMM.
The law $h$ is described graphically in figure~\ref{ISDA:fig:7}.
It is noteworthy that this model is equivalent to a simple HMM when $\Lambda=2$\,.
And when $\Lambda=1$\,, the law $h$ just maps the immediate observation to action, without any memory of the past observations.
\vspace{5pt}\\
For any $h\in\mathcal{H}$, define $P[h]$ the complete probabilistic law of the system Universe/Planner:
$$
\begin{array}{@{}l@{}}\displaystyle
P[h](x,y,z,m)=P(y_{1:T},z_{1:T}|x_{1:T})\prod_{t=1}^T h^0(x_t|m^1_{t})
\\\displaystyle\hspace{150pt}
\times h^1(m_t^1|y_{t-1},m^2_{t})\prod_{\lambda=2}^{\Lambda}h^\lambda(m_{t}^\lambda|m_{t-1}^{\lambda-1},m_{t}^{\lambda+1})\;.
\end{array}
$$
Our \emph{approximated} planning problem reduces to find the near-optimal strategy $h_O\in\mathcal{H}$ such that:
$$
%\begin{array}{@{}l@{}}\displaystyle
h_O\in\arg\max_{h\in\mathcal{H}}\sum_{x_{1:T}}\sum_{y_{1:T}}\sum_{z_{1:T}}\sum_{m^{1:\Lambda}_{1:T}}P[h](x,y,z,m)\,
%\\\displaystyle
%\hspace{140pt}\times 
V(x_{1:T},y_{1:T},z_{1:T})\;.
%\end{array}
$$
A solution to this problem, by means of the cross-entropy method, is proposed in the next section.
\section{Cross-entropic optimization of $h$}
The reader interested in CE methods should refer to the \emph{tutorial on the CE method}\cite{IEEEisda:boer}.
CE algorithms were first dedicated to estimating the probability of rare events.
A slight change of the basic algorithm made it also good for optimization.
In their new article\cite{MCO:mello}, Homem-de-Mello and Rubinstein have given some results about the global convergence.
In order to ensure such convergence, some refinements are introduced particularly about the selective rate.
\vspace{5pt}\\
This presentation is restricted to the CE optimization method.
The new improvements of the CE algorithm proposed in~\cite{MCO:mello} have not been implemented, but the algorithm has been seen to work properly.
For this reason, this paper does not deal with the choice of the selective rate.
\subsection{General CE algorithm for the optimization}
The Cross Entropy algorithm repeats until convergence the three successive phases:
\begin{enumerate}
\item Generate samples of random data according to a parameterized random mechanism,
\item Select the best samples according to an evaluation criterion,
\item Update the parameters of the random mechanism, on the basis of the selected samples.
\end{enumerate}
In the particular case of CE, the update in phase 3 is obtained by minimizing the Kullback-Leibler distance, or cross entropy, between the updated random mechanism and the selected samples.
The next paragraphs describe on a theoretical example how such method can be used in an optimization problem.
\paragraph{Formalism.}
Let be given a function $x\mapsto f(x)$; this function is easily computable.
The value $f(x)$ has to be maximized, by optimizing the choice of $x\in X$.
The function $f$ will be the evaluation criterion.
\vspace{5pt}\\
Now let be given a family of probabilistic laws, $P_\sigma|_{\sigma\in\Sigma}$\,, applying on the variable $x$.
The family $P$ is the parameterized random mechanism.
The variable $x$ is the random data.
\vspace{5pt}\\
Let $\rho\in\,]0,1[$ be a selective rate.
The CE algorithm for $(x,f,P)$ follows the synopsis\,:
\begin{enumerate}
\item Initialize $\sigma\in\Sigma$\,,
\item \label{XX:step2}Generate $N$ samples $x_n$ according to $P_\sigma$\,,
\item Select the $\rho N$ best samples according to the evaluation criterion $f$\,,
\item Update $\sigma$ as a minimizer of the cross-entropy with the selected samples:
$$\sigma\in\arg\max_{\sigma\in\Sigma}\sum_{n~{\rm selected}}\ln P_\sigma(x_n)\;,$$
\item Repeat from step \ref{XX:step2} until convergence.
\end{enumerate}
\emph{This algorithm requires $f$ to be easily computable.}
%\vspace{5pt}\\
\paragraph{Interpretation.}
The CE algorithm tightens the law $P_\sigma$ around the maximizer of $f$.
Then, when the probabilistic family $P$ is well suited to the maximization of $f$\,, it becomes equivalent to find a maximizer for $f$ or to optimize the parameter $\sigma$ by means of the CE algorithm.
The problem is to find a good family\dots
Another issue is the criterion for deciding the convergence.
Some answers are given in \cite{MCO:mello}.
Now, it is outside the scope of this paper to investigate these questions precisely.
Our criterion was to stop after a given threshold of successive \emph{unsuccessful tries} and this very simple method have worked fine on our problem.
\subsection{Application}\label{MCO:Subsect:MainAlgo}
Optimizing $h\in\mathcal{H}$ means tuning the parameter $h$ in order to tighten the probability $P[h]$ around the optimal values for $V$\,.
This is exactly solved by the \emph{Cross-Entropy} optimization method.
However, it is required that the evaluation function $V$ is easily computable.
Typically, the definition of $V$ may be recursive, \emph{eg.}\,:
$$
V(x_{1:T},y_{1:T},z_{1:T})=\left\{v_t(x_t,y_t,z_t,\right\}_{t=T}^{2}v_1(x_1,y_1,z_1)\left\{)\right\}_{t=T}^{2}\;.
$$
Let the \emph{selective rate} $\rho$ be a positive number such that $\rho<1$\,.
The cross-entropy method for optimizing $h$ follows the synopsis\,:
\begin{enumerate}
\item Initialize $h$\,. For example a flat $h$,
\item \label{Algo:0:1:FirstStep}Make $N$ tossing $\theta^n=(x^n,y^n,z^n,m^n)$ according to the law $P[h]$,
\item Choose the $\rho N$ best samples $\theta^n$ according to the sample evaluation $V(x^n_{1:T},y^n_{1:T},z^n_{1:T})$\,.
Denote $S$ the set of the selected samples,
\item Update $h$ as the minimizer of the cross-entropy with the selected samples:
\begin{equation}
\label{StratApp:1:1}
h\in\arg\max_{h\in\mathcal{H}}\sum_{n\in S}\ln P[h](\theta^n)\;,
\end{equation}
\item Reiterate from step~\ref{Algo:0:1:FirstStep} until convergence.
\end{enumerate}
For our HHMM model, the maximization~(\ref{StratApp:1:1}) is solved by:
$$
h^0(A|B)=\frac{
\mathrm{card}\Bigl\{n\in S\,,t\,/\,A=x_t^{n}\;~\mbox{and}~
B=m_t^{1,n}\Bigr\}
}{
\mathrm{card}\Bigl\{n\in S\,,t\,/\,B=m_t^{1,n}\Bigr\}
}\;,
$$
$$
h^1(A|B,C)=\frac{
\mathrm{card}\Bigl\{n\in S\,,t\,/\,A=m_t^{1,n}\,,
B=y_{t-1}^{n}\;~\mbox{and}~C=m_t^{2,n}\Bigr\}
}{
\mathrm{card}\Bigl\{n\in S\,,t\,/\,B=y_{t-1}^{n}\;~\mbox{and}~C=m_t^{2,n}\Bigr\}
}\;.
$$
and for $2\le\lambda\le\Lambda$\,,:
$$
h^\lambda(A|B,C)=\frac{
\mathrm{card}\Bigl\{n\in S\,,t\,/\,A=m_t^{\lambda,n}\,,
B=m_{t-1}^{\lambda-1,n}\;~\mbox{and}~C=m_t^{\lambda+1,n}\Bigr\}
}{
\mathrm{card}\Bigl\{n\in S\,,t\,/\,B=m_{t-1}^{\lambda-1,n}\;~\mbox{and}~C=m_t^{\lambda+1,n}\Bigr\}
}\;.
$$
\subsection{Bypassing the prior modelling}
POMDPs have a main drawback, they require the priors about the world and the mission evaluation to be defined.
Defining the mission evaluation is not a difficult task in general.
But it is quite uneasy to define the HMM which is modelling the world.
For example, what is the structure (memory requirement, time/space hierarchies) of this HMM?
What are the transition laws?
In many works, the structure is defined by the designer of the system.
But generally, the transition parameters are specified by learning from the actual world (\emph{eg.} Baum-Welch algorithm).
At last, the two step are implemented:
\begin{itemize}
\item Learn priorly the parameters of the HMM associated to the actual world (Baum-Welch),
\item Solve the POMDP for the instant play.
\end{itemize}
This is what is generally done when using the dynamic programming viewpoint.
But there is something nice with the cross-entropic methods.
Since the algorithm just needs samples, the \emph{HMM is not required}.
More precisely, since $P[h]=P(y_{1:T},z_{1:T}|x_{1:T})h(x_{1:T},m_{1:T}|y_{1:T})$\,, the update property~(\ref{StratApp:1:1}) may be factorized and the following simple result is derived:
$$
h\in\arg\max_{h\in\mathcal{H}}\sum_{n\in S}\ln h(x^n_{1:T},m^n_{1:T}|y^n_{1:T})\;,
$$
What is needed is a stochatic process, which, in association with the policy $h$\,, will build the samples.
This process may just be the \emph{actual world} itself!
Then, optimizing $h$ only requires time for experimental implementations\,:
\begin{enumerate}
\item Initialize $h$\,,
\item \label{Algo:X:0:1:FirstStep} $N$ samples $\theta^n=(x^n,y^n,m^n)$ and their respective evaluations $V^n$ are obtained by repeating $N$ times the following procedure:
\begin{enumerate}
\item Initialize $t=1$\,,
\item \label{Algo:X:0:1:StepTwo}Make a toss of $x^n_t$ and $m^n_t$ according to the law $h(x_t,m_t|m^n_{t-1},y^n_{t-1})$\,,
\item \emph{Execute the action $x^n_t$ in the actual world,}
\item \emph{Make the actual measurement $y^n_t$,}
\item Set $t=t+1$,
\item Repeat from step~\ref{Algo:X:0:1:StepTwo} until $t>T$,
\item Make an actual evaluation, denoted $V_n$, of the whole play (this evaluation depends on the actions, observations, and of the \emph{actual hidden states}),
\end{enumerate}
\item Choose the $\rho N$ best samples $\theta^n$ according to the sample evaluation $V^n$.
Denote $S$ the set of the selected samples,
\item Update $h$ as the minimizer of the cross-entropy with the selected samples:
$$
h\in\arg\max_{h\in\mathcal{H}}\sum_{n\in S}\ln h(x^n_{1:T},m^n_{1:T}|y^n_{1:T})\;,
$$
\item Reiterate from step~\ref{Algo:X:0:1:FirstStep} until convergence.
\end{enumerate}
This method avoids any construction of the world prior $P(y,z|x)$, but it requires many experimentations on the actual world.
This is sometimes unworkable.
\vspace{5pt}\\
The next section presents an example of implementation of the algorithm described in section~\ref{MCO:Subsect:MainAlgo}.
\section{Implementation}
\label{Isda:example}
The algorithm has been applied to a simplified target detection problem.
At this time, the \emph{open-loop} implementation has not been investigated yet.
\subsection{Problem setting}
A target $R$ is moving in a lattice of $20\times 20$ cells, \emph{ie.} $[\![\,0,19\,]\!]^2$.
$R$ is tracked by two mobiles, $B$ and $C$, controlled by the planner.
The coordinate of $R$, $B$ and $C$ at time $t$ are respectively denoted $(i_R^t,j_R^t),$ $(i_B^t,j_B^t)$ and $(i_C^t,j_C^t)$.
$B$ and $C$ have a very limited information about the target position, and are maneuvering much slower:\\
$\bullet$ A move for $B$ (respectively $C$) is either: \emph{turn left}, \emph{turn right}, \emph{go forward}, \emph{no move}.
Consequently, there are $4\times4=16$ possible actions for the planner.
These moves cannot be combined in a single turn. No diagonal forward: a mobile is either directed up, right, down or left,\\
$\bullet$ The mobiles are initially positioned in the down corners, \emph{ie.} $i^1_B=0,$ $j^1_B=19$ and $i^1_C=19,$ $j^1_C=19$.
The mobile are initially directed \emph{downward},\\
$\bullet$ $B$ (respectively $C$) observes whether the target relative position is forward or not.
More precisely:\vspace{2pt}\\
$\rien\quad\circ$ when $B$ is directed upward, the mobile knows whether $j_R<j_B$ or not,\\
$\rien\quad\circ$ when $B$ is directed right, the mobile knows whether $i_R>i_B$ or not,\\
$\rien\quad\circ$ when $B$ is directed downward, the mobile knows whether $j_R>j_B$ or not,\\
$\rien\quad\circ$ when $B$ is directed left, the mobile knows whether $i_R<i_B$ or not,\vspace{4pt}\\
$\bullet$ $B$ (respectively $C$) knows whether its distance with the target is less than $3$, \emph{ie.} $d_\infty(B,R)<3$, or not.
The distance $d_\infty$ is defined by:
$$
d_\infty(B,R)=\max\{|i_B-i_R|\,,|j_B-j_R|\}\;.
$$
At last, there are $2^4=16$ possible observations for the planner.
\vspace{5pt}\\
Several test cases have been considered.
In case 1, the target $R$ does not move.
In any other case, the target $R$ chooses stochastically its next position in its neighborhood.
Any move is possible (up/down, left/right, diagonals, no move).
The probability to choose a new position is proportional to the sum of the squared distance from the mobiles:
$$
\left\{
\begin{array}{@{}l@{}}\displaystyle
P(R^{t+1}|R^t)=0\mbox{ if }|i_{R}^{t+1}-i_{R}^{t}|>1\mbox{ or }|j_{R}^{t+1}-j_{R}^{t}|>1\;, \vspace{5pt}
\\\displaystyle
P(R_{t+1}|R_t)\propto(i_{R}^{t+1}-i_{B}^{t})^2+(j_{R}^{t+1}-j_{B}^{t})^2
\\\displaystyle\hspace{170pt}
+(i_{R}^{t+1}-i_{C}^{t})^2+(j_{R}^{t+1}-j_{C}^{t})^2\quad\mbox{else}\;.
\end{array}\right.
$$ 
This definition was intended to favorize escape moves: more great is a distance, more probable is the move.
But in such summation, a short distance will be neglected compared to a long distance (this was initially a mistake in the modelling!).
It is implied that a distant mobile will hide a nearby mobile.
This ``deluding'' property is interesting and will induce actually two different kinds of strategy, whithin the learned machines.
\vspace{5pt}\\
The objective of the planner is to maintain the target sufficiently closed to at least one mobile (in this example, the distance between the target and a mobile is required to be not more than $3$).
More precisely, the evaluation function, $V$, is just counting the number of such ``encounter'':
$$
V_0=0\ ;\quad V_{t}=V_{t-1}+1\mbox{\ \ if }d(B^t,R^t)\le3\mbox{ or }d(C^t,R^t)\le3\ ;\quad
V_{t}=V_{t-1}\mbox{\ \ else.}
$$
The total number of turns is $T=100$.
\subsection{Results}
\paragraph{Generality.}
Like many stochastic algorithms, this algorithm needs some time for convergence.
For the considered example, about two hours were needed for convergence (on a 2GHz PC); the selective rate was $\rho=0.5$.
This speed depends on the size of the HHMM model and on the convergence criterion.
A weak and a strong criterion are used for deciding the convergence.
Within the weak criterion, the algorithm is halted after $100$ successive unsuccessful tries.
Within the strong criterion, the algorithm is halted after $500$ successive unsuccessful tries.
Of course, the strong criterion computes a (slightly) better optimum than the weak criterion, but it needs time.
Because of the many tested examples, the weak criterion has been the most used in particular for the big models.
For the same HHMM model, the computed optimal values do not depend on the algorithmic instance (small variations result however from the stochastic nature of the algorithm):
the convergence seems to be global.
\vspace{5pt}\\
In the sequel, mean evaluations are rounded to the nearest integer.
The presentation is made clearer.
And owing to the small variations of this stochastic algorithm, more precision turns out to be irrelevant.
\paragraph{Case 1: $R$ does not move.}
This example has been considered in order to test the algorithm.
The position of the target is fixed in the center of the square space, \emph{ie.} $i^1_R=j^1_R=10$.
It is recalled that the mobiles are initially directed downward.
Then, the optimal strategy is known and its value is $85$\,: the time needed to reach the target is $15$\,, and no further move is needed.
The learned $h_O$ approximates the evaluation $84$\,.
The convergence is good.
\paragraph{Case 2: $R$ is moving but the observation $y$ is hidden.}
Initially, $R$ is located within the $20\times10$ upper cells of the lattice (\emph{ie.} $[\![\,0,19\,]\!]\times[\![\,0,9\,]\!]$), accordingly to a uniform probabilistic law.
The computed optimal means evaluation is about $32$.
In this case, the mobiles tend to move towards the upper corners.
\begin{figure}
\caption{Near-optimal control sequence}
\label{ISDA:fig:8}
\begin{center}
\begin{tabular}{l}\vspace{-15pt}\\
\fbox{\scalebox{.8}{\xymatrix{
\bullet^{5,6}\ar[r]&\bullet^{0,7}\ar[d]\times^0\ar[r]&\bullet^8\ar[r]\times^1\ar[dr]&\bullet^9\ar[r]&\bullet^{10}\times^3\ar[dr]&&&&\\
\bullet^{3,4}\ar[u]&\bullet^{1,2}\ar[l]&&\times^2\ar[ur]&&\times^4\ar[r]&\times^{5,6}\ar[d]&&\\
&&&&&&\times^7\ar[r]&\times^8\ar[d]&\times^{10}\\
&&&&&&&\times^9\ar[ru]&\\
&&&&&&&\circ^{10}&\\
&&\circ^0\ar[d]&&&&&\circ^9\ar[u]&\\
&&\circ^{1,2}\ar[r]&\circ^{3}\ar[r]&\circ^{4}\ar[r]&\circ^{5}\ar[r]&\circ^{6}\ar[r]&\circ^{7,8}\ar[u]&
}}}\\\\
%\begin{tabular}{l}
$\times=\;$target\qquad
$\bullet=\;$observer 1\qquad
$\circ=\;$observer 2
\\
Relative times are put in supscript\vspace{-10pt}
\end{tabular}
\end{center}
\end{figure}
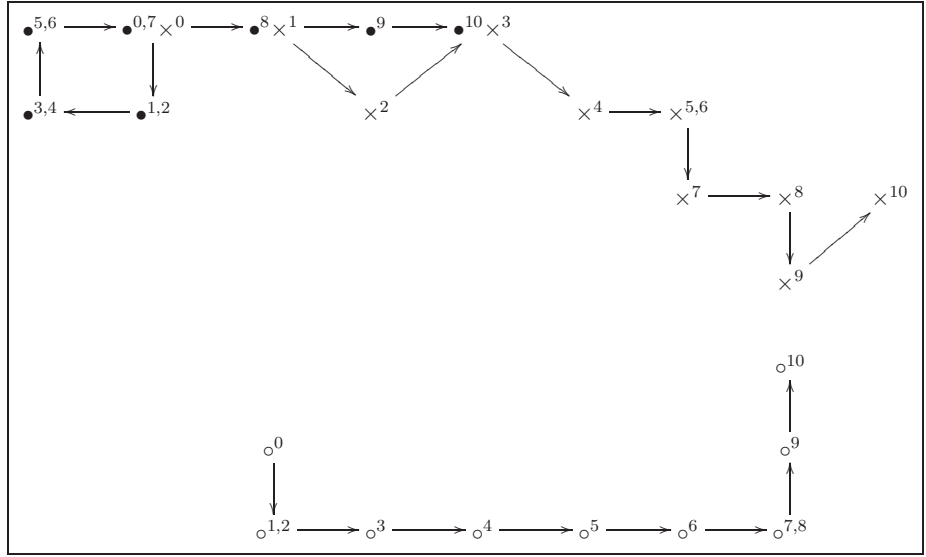
\paragraph{Case 3: $R$ is moving and $y$ is observed.}
Again, $R^1$ is located uniformly within the $20\times10$ upper cells of the lattice.
The computed optimal means evaluation is about $69$.
This evaluation has been obtained from a large HHMM model ($\Lambda=2$ with $256$ states per level, \emph{ie.} $\mathrm{card}(M^\lambda)=256$) and with the strong criterion.
However, somewhat smaller models should work as well.
\vspace{3pt}\\
Specific computations are now presented, depending on the number of levels $\Lambda$ and the number of states per levels.
For each case, the weak criterion has been used.
\vspace{5pt}\\
\emph{Subcase $\Lambda=1$\,.}
For such model, the action $x_t$ is constructed only from the immediate last observation $y_{t-1}$.
The model does not keep any memory of the past observations.
Then, only $16$ states are sufficient to describe the hidden variable $m_t^1$\,, \emph{ie.} $\mathrm{card}(M^1)=16$.
The resulting optimum is $54$\,.
\vspace{5pt}\\
\emph{Subcase $\Lambda=2$\,.}
This model is equivalent to a HMM and it is assumed that $\mathrm{card}(M^1)=\mathrm{card}(M^2)$.
The following table gives the computed optimum for several choices of the set of states:
$$
\begin{array}{|l||l|l|l|l|}
\hline \mathrm{card}(M^\lambda)&16&32&64&256
\\
\hline\mathrm{Evaluation}&65&66&67&67
\\\hline
\end{array}
$$
%\newpage\rien\vspace{200pt}\\
%
%\newpage
It is noteworthy that the memory of the past observations allows better strategies than just immediate observations (case $\Lambda=1$)\,.
Indeed, the evaluation jumps from $54$ up to $67$.
\vspace{5pt}\\
\emph{Subcases $\Lambda>2$.}
A comparison of graduated hierarchic models, $1\le\Lambda\le4$, has been made.
The first level contained $16$ possible states, and the higger levels were restricted to $2$ states:
$$
\begin{array}{|l||l|l|l|l|}
\hline\mathrm{hierarchic\ grade}&\Lambda=1&\Lambda=2&\Lambda=3&\Lambda=4
\\
\hline\mathrm{card}(M^\lambda)|_{\lambda=1}^\Lambda&16&16,2&16,2,2&16,2,2,2
\\\hline
\end{array}
$$
The test has been accomplished according to the weak criterion:
$$
\begin{array}{|l||l|l|l|l|}
\hline\mathrm{hierarchic\ grade}&\Lambda=1&\Lambda=2&\Lambda=3&\Lambda=4
\\
\hline\mathrm{Evaluation\ (weak)}&54&59&\bf 56&65
\\\hline
\end{array}
$$
and the strong criterion:
$$
\begin{array}{|l||l|l|l|l|}
\hline\mathrm{hierarchic\ grade}&\Lambda=1&\Lambda=2&\Lambda=3&\Lambda=4
\\
\hline\mathrm{Evaluation\ (strong)}&55&61&64&66
\\\hline
\end{array}
$$
It seems that a high hierarchic grade makes the convergence more difficult.
This is particularly the case here for the grade $\Lambda=3$\,, which failed under the weak criterion at only $56$.
However, the algorithm still works well by improving the convergence criterion.
\vspace{5pt}\\
It is interesting to make a comparison with the subcase $\Lambda=2$ where $\mathrm{card}(M^1)=\mathrm{card}(M^2)=16$.
Under the weak criterion, the result for this HHMM was $65$ as for the grade $\Lambda=4$.
However, the dimension of the law is quite different for the two models:
\begin{itemize}
\item $15\times 16+15\times 16\times 16+15\times 16=4320$ for the $2$-level HHMM,
\item $15\times 16+15\times 16\times 2+1\times 16\times 2+1\times 2\times 2+1\times 2=758$ for the $4$-level HHMM.
\end{itemize}
This dimension is a rough characterization of the complexity of the model.
It seems clear on these examples that the highly hierarchized models are more efficient than the weakly hierarchized models.
And the problem considered here is quite simple.
On complex problems, hierarchical models may be pre-eminent.
\paragraph{Global behavior.}
\rien\\\emph{The algorithm.}
The convergence speed is low at the beginning.
After this initial stage, it improves greatly until it reaches a new ``waiting'' stage.
This alternation of low speed and great speed stages have been noticed several times until an acceptable convergence.
Nevertheless, the speed is globally decreasing with the time.
\vspace{5pt}\\
\emph{The near optimal policy.}
It is now discussed about the behaviour of the best found policy.
This HHMM has reach the mean evaluation $69$.
The mobiles strategy results in a tracking of the target.
The figure~\ref{ISDA:fig:8} illustrates a short sequence of escape/tracking of the target.
It has been noticed two quite distinct behaviours, among the many runs of the policy:
\begin{itemize}
\item \emph{The two mobiles may both cooperate on tracking the target,}
\item \emph{When the target is near the border, one mobile may stay along the opposite border while the other mobile may perform the tracking.}
This strategy seems strange at first sight.
But it is recalled that the moving rule of the target tends to neglect a nearby mobile compared to a distant mobile.
In this strategy, the first mobile is just annihilating the ability of the target to escape from the tracking of the second mobile.
\end{itemize}
\section{Conclusion}
In this paper, we proposed a general method for approximating the optimal planning in a partially observable control problem.
Hierarchical HMM families have been used for approximating the optimal decision tree, and the approximation has been optimized by means of the Cross-Entropy method.
Moreover, hierarchical HMM has been characterized as a multilevel HMM with alternately up-and-down (Markovian) diffusions of the information between the levels.
The algorithm has been applied on a simplier model based on a weaker diffusion of the information.
\vspace{5pt}\\
At this time, the method has been applied to a strictly discrete-state problem and has been seen to work properly.
The convergence seems global: the many runs of the algorithm have reach the approximately same optimal values.
The optimal HHMMs are able to ``track'' the target.
An interesting point is that these HHMMs have discovered two quite different global strategies and are able to choose between them: make the mobiles both cooperate on tracking or require one mobile for deluding the target.
\vspace{5pt}\\
The results are promising, but the example test is still simple.
The observation and action spaces are limited to a few number of states.
And what happens if the hidden space becomes much more intricated?
There are several possible answers to such difficulties:
\vspace{5pt}\\
First, it has been shown that it is not necessary to define the priors of the problem (defining priors often implies many approximations):
the CE algorithm is able to learn the policy directly from the actual world.
Moreover, the cross-entropic principle could be applied for optimizing continuous laws.
It is thus certainly possible to consider mixed continuous/discrete HHMM, which are more realistic for a planning policy.
At last, many refinements are foreseeable about the structure of the HHMMs.
More precisely, hierarchic models for observation, action and memory (hierarchical HMMs, but also hierarchical Markov Random Fields) should be improved in order to locally factorize intricated problems.
This research is just preliminary and future works should investigate these questions.
\appendix
\section{Murphy and Paskin Bayesian Network}
\label{sect:App:hhmm}
\begin{figure}
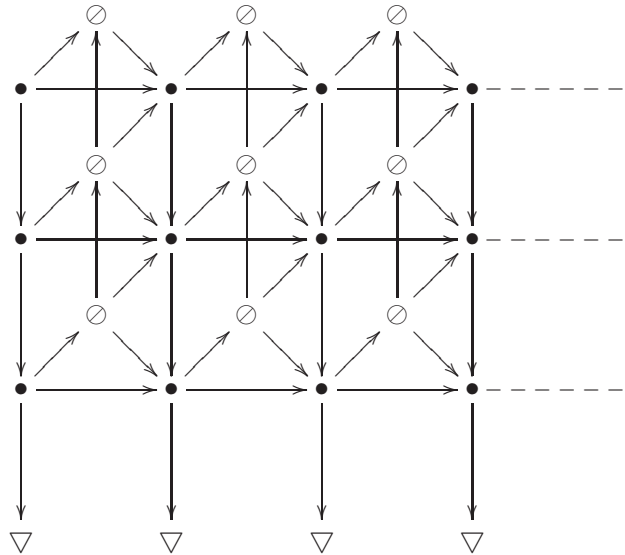

\caption{A Hierarchical HMM (4 levels)}
\label{figMFcp:6}
\begin{center}\begin{tabular}{@{}c@{}}
\rien\vspace{-45pt}\\
\xy
(10,-10)*+{\oslash}="s00",
(30,-10)*+{\oslash}="s01",
(50,-10)*+{\oslash}="s02",
(0,-20)*+{\bullet}="b10",
(20,-20)*+{\bullet}="b11",
(40,-20)*+{\bullet}="b12",
(60,-20)*+{\bullet}="b13",
(80,-20)="bf1",
\ar @{->}"b10";"s00",
\ar @{<-}"b11";"s00",
\ar @{->}"b11";"s01",
\ar @{<-}"b12";"s01",
\ar @{->}"b12";"s02",
\ar @{<-}"b13";"s02",
\ar @{->}"b10";"b11",
\ar @{->}"b11";"b12",
\ar @{->}"b12";"b13",
\ar @{--}"b13";"bf1",\POS "bf1",
(10,-30)*+{\oslash}="s10",
(30,-30)*+{\oslash}="s11",
(50,-30)*+{\oslash}="s12",
\ar @{<-}"b11";"s10",
\ar @{<-}"b12";"s11",
\ar @{<-}"b13";"s12",
\ar @{->}"s10";"s00",
\ar @{->}"s11";"s01",
\ar @{->}"s12";"s02",
\POS "s02",
(0,-40)*+{\bullet}="b00",
(20,-40)*+{\bullet}="b01",
(40,-40)*+{\bullet}="b02",
(60,-40)*+{\bullet}="b03",
(80,-40)="bf0",
\ar @{->}"b10";"b00",
\ar @{->}"b11";"b01",
\ar @{->}"b12";"b02",
\ar @{->}"b13";"b03",
\ar @{->}"b00";"s10",
\ar @{<-}"b01";"s10",
\ar @{->}"b01";"s11",
\ar @{<-}"b02";"s11",
\ar @{->}"b02";"s12",
\ar @{<-}"b03";"s12",
\ar @{->}"b00";"b01",
\ar @{->}"b01";"b02",
\ar @{->}"b02";"b03",
\ar @{--}"b03";"bf0",\POS "bf0",
(10,-50)*+{\oslash}="s00",
(30,-50)*+{\oslash}="s01",
(50,-50)*+{\oslash}="s02",
\ar @{<-}"b01";"s00",
\ar @{<-}"b02";"s01",
\ar @{<-}"b03";"s02",
\ar @{<-}"s10";"s00",
\ar @{<-}"s11";"s01",
\ar @{<-}"s12";"s02",
\POS "s02",
(0,-60)*+{\bullet}="b10",
(20,-60)*+{\bullet}="b11",
(40,-60)*+{\bullet}="b12",
(60,-60)*+{\bullet}="b13",
(80,-60)="bf1",
\ar @{<-}"b10";"b00",
\ar @{<-}"b11";"b01",
\ar @{<-}"b12";"b02",
\ar @{<-}"b13";"b03",
\ar @{->}"b10";"s00",
\ar @{<-}"b11";"s00",
\ar @{->}"b11";"s01",
\ar @{<-}"b12";"s01",
\ar @{->}"b12";"s02",
\ar @{<-}"b13";"s02",
\ar @{->}"b10";"b11",
\ar @{->}"b11";"b12",
\ar @{->}"b12";"b13",
\ar @{--}"b13";"bf1",\POS "bf1",
(0,-80)*+{\bigtriangledown}="b00",
(20,-80)*+{\bigtriangledown}="b01",
(40,-80)*+{\bigtriangledown}="b02",
(60,-80)*+{\bigtriangledown}="b03",
\ar @{->}"b10";"b00",
\ar @{->}"b11";"b01",
\ar @{->}"b12";"b02",
\ar @{->}"b13";"b03", \POS "b03",
\endxy\vspace{-10pt}
\end{tabular}
\end{center}
\end{figure}

\begin{figure}
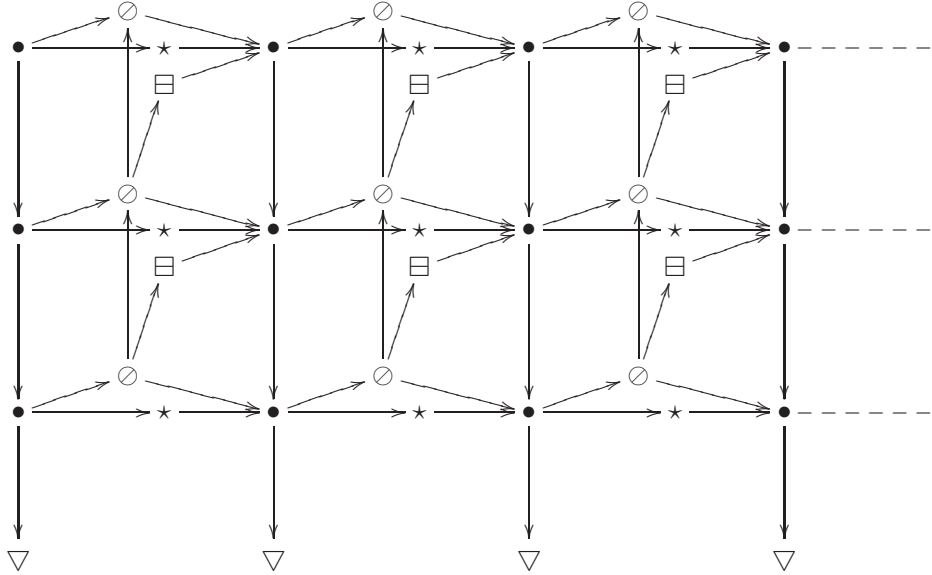

\caption{Simplified Model of a HHMM}
\label{MCO:fig:11}
\begin{center}
\begin{tabular}{@{}c@{}}
\vspace{-60pt}\\
\xy
(-.5,-14.5)*+{},
%\green
(0,-15)*+{\bullet}="b10",
(15,-15)*+{\circ}="b11",
(30,-15)*+{\bullet}="b12",
(45,-15)*+{\circ}="b13",
(60,-15)*+{\bullet}="b14",
(75,-15)="b15",
\ar @{->}"b10";"b11",
\ar @{->}"b11";"b12",
\ar @{->}"b12";"b13",
\ar @{->}"b13";"b14",
\ar @{--}"b14";"b15",\POS "b15",
(0,-30)*+{\bullet}="b00",
(15,-30)*+{\circ}="b01",
(30,-30)*+{\bullet}="b02",
(45,-30)*+{\circ}="b03",
(60,-30)*+{\bullet}="b04",
(75,-30)="b05",
\ar @{<-}"b00";"b10",
\ar @{->}"b01";"b11",
\ar @{<-}"b02";"b12",
\ar @{->}"b03";"b13",
\ar @{<-}"b04";"b14",
\ar @{->}"b00";"b01",
\ar @{->}"b01";"b02",
\ar @{->}"b02";"b03",
\ar @{->}"b03";"b04",
\ar @{--}"b04";"b05",\POS "b05",
(0,-45)*+{\bullet}="b10",
(15,-45)*+{\circ}="b11",
(30,-45)*+{\bullet}="b12",
(45,-45)*+{\circ}="b13",
(60,-45)*+{\bullet}="b14",
(75,-45)="b15",
\ar @{->}"b00";"b10",
\ar @{<-}"b01";"b11",
\ar @{->}"b02";"b12",
\ar @{<-}"b03";"b13",
\ar @{->}"b04";"b14",
\ar @{->}"b10";"b11",
\ar @{->}"b11";"b12",
\ar @{->}"b12";"b13",
\ar @{->}"b13";"b14",
\ar @{--}"b14";"b15",\POS "b15",
(0,-60)*+{\bigtriangledown}="b00",
(30,-60)*+{\bigtriangledown}="b02",
(60,-60)*+{\bigtriangledown}="b04",
\ar @{<-}"b00";"b10",
\ar @{<-}"b02";"b12",
\ar @{<-}"b04";"b14",\POS "b14",
(75.5,-60.5)*+{}
\endxy\vspace{10pt}\\
$\circ$ contains $\star$ , $\oslash$ and possibly $\boxminus$\vspace{-10pt}
\end{tabular}
\end{center}
\end{figure}

\paragraph{Definition.} A HHMM is a HMM which output is either an observation (when reaching the output level) or a HHMM.
A formal definition is given here, for peoples acquainted with such HMM formalism.
But it is not explained in detail.
\subparagraph{Formalism.} A HHMM with $\Delta$ levels is characterized by an observation set $Q_1$ and by $\Delta-1$ quadruplet $(Q_d,E_d,\rho_d,\pi_d)_{2\le d\le \Delta}$.
For each level \mbox{$d\in[\![2,\Delta]\!]$}\,, the quadruplet $(Q_d,E_d,\rho_d,\pi_d)$ is related to a HMM for the level $d$, thus verifying:
\begin{itemize}
\item $Q_d$ is a set of states for the level $d$,
\item $E_d\subset Q_d$ is a set of ending states (may be empty for the root $\Delta$).
When such state is obtained, the HMM stops and jumps back to the HMM of higher level; for the root level, the HHMM just ends,
\item $\rho_d(q_{d-1,t},q_{d,t})=P(q_{d-1,t}|q_{d,t})$ describes the probability for the HMM at level $d$ to produce (at time $t$) the output $q_{d-1,t}\in Q_{d-1}$ when the inner state is $q_{d,t}\in Q_d$.
When $d=2$, the output is just an observation produced by the HHMM.
When $d>2$, the HMM of level $d$ initiates the children HMM of level $d-1$ with the starting state $q_{d-1,t}$\,, and waits until the children HMM reachs an ending state,
\item $\pi_d(q_{d,t},q_{d,t+1})=P(q_{d,t+1}|q_{d,t})$ describes the probability for the HMM at level $d$ to transit from a state $q_{d,t}\in Q_d\setminus E_d$ at time $t$ to the state $q_{d,t+1}\in Q_d$ at time $t+1$.
This transit is runned \emph{after} the production stage.
\end{itemize}
In this definition, $\pi_\Delta(\emptyset,q_{\Delta,1})=P(q_{\Delta,1})$ initializes the HHMM at time $1$.
\paragraph{Related Bayesian Network.}
Murphin and Paskin\cite{IEEEisda:murphy} proposed an alternate definition of the HHMM by means of Bayesian Network.
This BN is in fact a vectorized HMM, which has the Markov property both in time and in hierarchical level.
This BN is relying on $2$ type of cells:
state cells (symbols $\bullet,\bigtriangledown$) and boolean cells (symbols $\oslash,\odot,\otimes$).
More precisely:
\begin{center}
\begin{tabular}{@{}l@{}}
$\bullet$ means an inner state cell,
\\
$\bigtriangledown$ means an output state cell,
\\
$\oslash$ means a \emph{unspecified} boolean cell,
\\
$\odot$ means a boolean cell \emph{specified} FALSE,
\\
$\otimes$ means a boolean cell \emph{specified} TRUE.
\end{tabular}
\end{center}
The structure of the Bayesian Network is given by figure~\ref{figMFcp:6}.
The boolean cells are dedicated to the control of the state transition.
A $\odot$/FALSE cell indicates that the bottom level is running.
A $\otimes$/TRUE cell indicates that the bottom level is ended.
Depending on the boolean configuration, it is chosen a HMM transit $\pi$\,, a HMM production $\rho$ or a wait (the states are leaved unchanged by the identity $\mathrm{id}$).
All the state transitions are summarized in tables~\ref{tabMFcp:2} and~\ref{tabMFcp:2:1}.
For example, when a children level is running, the father level is waiting and the father row is between two $\odot$ cells: the transition is $\mathrm{id}$ and the father state is leaved unchanged 
\vspace{5pt}\\
On the other hand, the table~\ref{tabMFcp:3} summarizes the boolean transitions.
A boolean cell may just be a copy of the bottom boolean cell, when this cell is False.
Indeed, that means that the children level is running, and consequently the father level is running although in a waiting mode.
Otherwise, the transition is controlled by the HMM transit $\pi$\,.
\begin{table}
\caption{Transition rule: for a state cell}
\label{tabMFcp:2}
\begin{center}
\begin{tabular}{l||c|c|c|c|}
Config.&\xy
(-5,0)*+{},
(5,0)*+{},
(0,10)*+{\bullet}="x01",
(0,-10)*+{\bullet/\bigtriangledown}="x21",
\ar @{->}"x01";"x21",
\POS "x21"
\endxy
&
\xy
%(-2.5,2.5)*+{},
(20,10)*+{\bullet}="x00",
(10,0)*+{\odot}="x10",
(0,-10)*+{\bullet}="x01",
(20,-10)*+{\bullet}="x21",
(10,-20)*+{\odot}="x12",
\ar @{->}"x00";"x21",
\ar @{->}"x10";"x21",
\ar @{->}"x01";"x21",
\ar @{<-}"x21";"x12",
\POS "x12",
%%
%(22.5,-22.5)*+{}
\endxy
&
\xy
%(-2.5,2.5)*+{},
(20,10)*+{\bullet}="x00",
(10,0)*+{\odot}="x10",
(0,-10)*+{\bullet}="x01",
(20,-10)*+{\bullet}="x21",
(10,-20)*+{\otimes}="x12",
\ar @{->}"x00";"x21",
\ar @{->}"x10";"x21",
\ar @{->}"x01";"x21",
\ar @{<-}"x21";"x12",
\POS "x12",
%%
%(22.5,-22.5)*+{}
\endxy
&
\xy
%(-2.5,2.5)*+{},
(20,10)*+{\bullet}="x00",
(10,0)*+{\otimes}="x10",
(0,-10)*+{\bullet}="x01",
(20,-10)*+{\bullet}="x21",
(10,-20)*+{\oslash}="x12",
\ar @{->}"x00";"x21",
\ar @{->}"x10";"x21",
\ar @{->}"x01";"x21",
\ar @{<-}"x21";"x12",
\POS "x12",
%%
%(22.5,-22.5)*+{}
\endxy
\\\vspace{-10pt}&&&&\\\hline\vspace{-10pt}&&&\\
Result.&
\xy
(0,10)*+{\bullet}="x01",
(0,-10)*+{\bullet/\bigtriangledown}="x21",
\ar^{\rho} @{->}"x01";"x21",
\POS "x21"
\endxy
&
\xy
%(-2.5,2.5)*+{},
(10,10)*+{}="x10",
(0,0)*+{\bullet}="x01",
(20,0)*+{\bullet}="x21",
(10,-10)*+{}="x12",
\ar^{\mathrm{id}} @{->}"x01";"x21",
\POS "x21",
%%
%(22.5,-22.5)*+{}
\endxy
&
\xy
%(-2.5,2.5)*+{},
(10,10)*+{}="x10",
(0,0)*+{\bullet}="x01",
(20,0)*+{\bullet}="x21",
(10,-10)*+{}="x12",
\ar^{\pi} @{->}"x01";"x21",
\POS "x21",
%%
%(22.5,-22.5)*+{}
\endxy
&
\xy
(0,10)*+{\bullet}="x01",
(0,-10)*+{\bullet}="x21",
\ar^{\rho} @{->}"x01";"x21",
\POS "x21"
\endxy
\end{tabular}
\end{center}
\end{table}
\begin{table}
\caption{Transition rule: for a state cell (border cases)}
\label{tabMFcp:2:1}
\begin{center}
\begin{tabular}{l||c|c|c|c|}
Config.&
\xy
%(-2.5,2.5)*+{},
(10,10)*+{\odot}="x10",
(0,0)*+{\bullet}="x01",
(20,0)*+{\bullet}="x21",
(10,-10)*+{\odot}="x12",
\ar @{->}"x10";"x21",
\ar @{->}"x01";"x21",
\ar @{<-}"x21";"x12",
\POS "x12",
%%
%(22.5,-22.5)*+{}
\endxy
&
\xy
%(-2.5,2.5)*+{},
(10,10)*+{\odot}="x10",
(0,0)*+{\bullet}="x01",
(20,0)*+{\bullet}="x21",
(10,-10)*+{\otimes}="x12",
\ar @{->}"x10";"x21",
\ar @{->}"x01";"x21",
\ar @{<-}"x21";"x12",
\POS "x12",
%%
%(22.5,-22.5)*+{}
\endxy
&
\xy
%(-2.5,2.5)*+{},
(20,10)*+{\bullet}="x00",
(10,0)*+{\odot}="x10",
(0,-10)*+{\bullet}="x01",
(20,-10)*+{\bullet}="x21",
\ar @{->}"x00";"x21",
\ar @{->}"x10";"x21",
\ar @{->}"x01";"x21",
\POS "x21",
%%
%(22.5,-22.5)*+{}
\endxy
&
\xy
%(-2.5,2.5)*+{},
(20,10)*+{\bullet}="x00",
(10,0)*+{\otimes}="x10",
(0,-10)*+{\bullet}="x01",
(20,-10)*+{\bullet}="x21",
\ar @{->}"x00";"x21",
\ar @{->}"x10";"x21",
\ar @{->}"x01";"x21",
\POS "x21",
%%
%(22.5,-22.5)*+{}
\endxy
\\\vspace{-10pt}&&&&\\\hline\vspace{-10pt}&&&\\
Result.&
\xy
%(-2.5,2.5)*+{},
(10,10)*+{}="x10",
(0,0)*+{\bullet}="x01",
(20,0)*+{\bullet}="x21",
(10,-10)*+{}="x12",
\ar^{\mathrm{id}} @{->}"x01";"x21",
\POS "x21",
%%
%(22.5,-22.5)*+{}
\endxy
&
\xy
%(-2.5,2.5)*+{},
(10,10)*+{}="x10",
(0,0)*+{\bullet}="x01",
(20,0)*+{\bullet}="x21",
(10,-10)*+{}="x12",
\ar^{\pi} @{->}"x01";"x21",
\POS "x21",
%%
%(22.5,-22.5)*+{}
\endxy
&
\xy
%(-2.5,2.5)*+{},
(10,10)*+{}="x10",
(0,0)*+{\bullet}="x01",
(20,0)*+{\bullet}="x21",
(10,-10)*+{}="x12",
\ar^{\pi} @{->}"x01";"x21",
\POS "x21",
%%
%(22.5,-22.5)*+{}
\endxy
&
\xy
(0,10)*+{\bullet}="x01",
(0,-10)*+{\bullet}="x21",
\ar^{\rho} @{->}"x01";"x21",
\POS "x21"
\endxy
\end{tabular}
\end{center}
\end{table}
\begin{table}
\caption{Transition rule: for a boolean cell}
\label{tabMFcp:3}
\begin{center}
\begin{tabular}{l||c|c|c|}
Config.&
\xy
%(-2.5,2.5)*+{},
(10,10)*+{\oslash}="x10",
(0,0)*+{\bullet}="x01",
(10,-10)*+{\odot}="x12",
\ar @{->}"x01";"x10",
\ar @{<-}"x10";"x12",
\POS "x12",
%%
%(22.5,-22.5)*+{}
\endxy
&
\xy
%(-2.5,2.5)*+{},
(10,10)*+{\oslash}="x10",
(0,0)*+{\bullet}="x01",
(10,-10)*+{\otimes}="x12",
\ar @{->}"x01";"x10",
\ar @{<-}"x10";"x12",
\POS "x12",
%%
%(22.5,-22.5)*+{}
\endxy
&
\xy
%(-2.5,2.5)*+{},
(10,10)*+{\oslash}="x10",
(0,0)*+{\bullet}="x01",
\ar @{->}"x01";"x10",
\POS "x10",
%%
%(22.5,-22.5)*+{}
\endxy
\\\vspace{-10pt}&&&\\\hline\vspace{-10pt}&&&\\
Result.&
$\oslash:=\odot$
&
\xy
%(-2.5,2.5)*+{},
(10,5)*+{\oslash}="x10",
(0,-5)*+{\bullet}="x01",
\ar^\pi @{->}"x01";"x10",
\POS "x10",
%%
%(22.5,-22.5)*+{}
\endxy
&
\xy
%(-2.5,2.5)*+{},
(10,5)*+{\oslash}="x10",
(0,-5)*+{\bullet}="x01",
\ar^\pi @{->}"x01";"x10",
\POS "x10",
%%
%(22.5,-22.5)*+{}
\endxy
\end{tabular}
\end{center}
\end{table}
\paragraph{Bayesian Network equivalence.}
It is shown now that the BN of figure~\ref{figMFcp:6} is a particular case of the BN of figure~\ref{ISDA:fig:6}\,.
The simple hierarchical model of figure~\ref{ISDA:fig:6} is thus sufficient to describes a general HHMM.
\vspace{5pt}\\
The BN of Murphy is really not far from figure~\ref{ISDA:fig:6}.
It is just needed to distort the arrows of figure~\ref{figMFcp:6}, to add intermediate state cells, denoted $\star$, and intermediate boolean cells, denoted $\boxminus$.
The intermediate Bayesian Network of figure~\ref{MCO:fig:10} is obtained then, and is clearly equivalent to the BN of Murphy.
By fusing the neighbooring cells $\oslash,\star,\boxminus$ into the cells $\circ$ and redefining the transition rules properly, the BN of figure~\ref{MCO:fig:11} is resulting, which is of the same kind as in figure~\ref{ISDA:fig:6}.
Notice that the cells $\circ$ are possibly containing two booleans and an information.
\begin{figure}
\caption{Modified Hierarchical HMM}
\label{MCO:fig:10}
\begin{center}
\begin{tabular}{@{}c@{}}
\rien\vspace{-60pt}\\
\scalebox{.97}{
\xy
(15,-15)*+{\oslash}="s00",
(50,-15)*+{\oslash}="s01",
(85,-15)*+{\oslash}="s02",
(0,-20)*+{\bullet}="b10",
(35,-20)*+{\bullet}="b11",
(70,-20)*+{\bullet}="b12",
(105,-20)*+{\bullet}="b13",
(125,-20)="bf1",
(20,-20)*+{\star}="ib0",
(55,-20)*+{\star}="ib1",
(90,-20)*+{\star}="ib2",
\ar @{->}"b10";"s00",
\ar @{<-}"b11";"s00",
\ar @{->}"b11";"s01",
\ar @{<-}"b12";"s01",
\ar @{->}"b12";"s02",
\ar @{<-}"b13";"s02",
\ar @{->}"b10";"ib0",
\ar @{->}"ib0";"b11",
\ar @{->}"b11";"ib1",
\ar @{->}"ib1";"b12",
\ar @{->}"b12";"ib2",
\ar @{->}"ib2";"b13",
\ar @{--}"b13";"bf1",\POS "bf1",
(15,-40)*+{\oslash}="s10",
(50,-40)*+{\oslash}="s11",
(85,-40)*+{\oslash}="s12",
(20,-25)*+{\boxminus}="is0",
(55,-25)*+{\boxminus}="is1",
(90,-25)*+{\boxminus}="is2",
\ar @{<-}"b11";"is0",
\ar @{<-}"is0";"s10",
\ar @{<-}"b12";"is1",
\ar @{<-}"is1";"s11",
\ar @{<-}"b13";"is2",
\ar @{<-}"is2";"s12",
\ar @{->}"s10";"s00",
\ar @{->}"s11";"s01",
\ar @{->}"s12";"s02",
\POS "s02",
(0,-45)*+{\bullet}="b00",
(35,-45)*+{\bullet}="b01",
(70,-45)*+{\bullet}="b02",
(105,-45)*+{\bullet}="b03",
(125,-45)="bf0",
(20,-45)*+{\star}="ib0",
(55,-45)*+{\star}="ib1",
(90,-45)*+{\star}="ib2",
\ar @{->}"b10";"b00",
\ar @{->}"b11";"b01",
\ar @{->}"b12";"b02",
\ar @{->}"b13";"b03",
\ar @{->}"b00";"s10",
\ar @{<-}"b01";"s10",
\ar @{->}"b01";"s11",
\ar @{<-}"b02";"s11",
\ar @{->}"b02";"s12",
\ar @{<-}"b03";"s12",
\ar @{->}"b00";"ib0",
\ar @{->}"ib0";"b01",
\ar @{->}"b01";"ib1",
\ar @{->}"ib1";"b02",
\ar @{->}"b02";"ib2",
\ar @{->}"ib2";"b03",
\ar @{--}"b03";"bf0",\POS "bf0",
(15,-65)*+{\oslash}="s00",
(50,-65)*+{\oslash}="s01",
(85,-65)*+{\oslash}="s02",
(20,-50)*+{\boxminus}="is0",
(55,-50)*+{\boxminus}="is1",
(90,-50)*+{\boxminus}="is2",
\ar @{<-}"b01";"is0",
\ar @{<-}"is0";"s00",
\ar @{<-}"b02";"is1",
\ar @{<-}"is1";"s01",
\ar @{<-}"b03";"is2",
\ar @{<-}"is2";"s02",
\ar @{<-}"s10";"s00",
\ar @{<-}"s11";"s01",
\ar @{<-}"s12";"s02",
\POS "s02",
(0,-70)*+{\bullet}="b10",
(35,-70)*+{\bullet}="b11",
(70,-70)*+{\bullet}="b12",
(105,-70)*+{\bullet}="b13",
(125,-70)="bf1",
(20,-70)*+{\star}="ib0",
(55,-70)*+{\star}="ib1",
(90,-70)*+{\star}="ib2",
\ar @{<-}"b10";"b00",
\ar @{<-}"b11";"b01",
\ar @{<-}"b12";"b02",
\ar @{<-}"b13";"b03",
\ar @{->}"b10";"s00",
\ar @{<-}"b11";"s00",
\ar @{->}"b11";"s01",
\ar @{<-}"b12";"s01",
\ar @{->}"b12";"s02",
\ar @{<-}"b13";"s02",
\ar @{->}"b10";"ib0",
\ar @{->}"ib0";"b11",
\ar @{->}"b11";"ib1",
\ar @{->}"ib1";"b12",
\ar @{->}"b12";"ib2",
\ar @{->}"ib2";"b13",
\ar @{--}"b13";"bf1",\POS "bf1",
(0,-90)*+{\bigtriangledown}="b00",
(35,-90)*+{\bigtriangledown}="b01",
(70,-90)*+{\bigtriangledown}="b02",
(105,-90)*+{\bigtriangledown}="b03",
\ar @{->}"b10";"b00",
\ar @{->}"b11";"b01",
\ar @{->}"b12";"b02",
\ar @{->}"b13";"b03", \POS "b03",
\endxy}
\end{tabular}
\end{center}
\end{figure}
%\section{Global convergence of the CE algorithm}
%\label{sect:App:CE}
%
% that's all folks
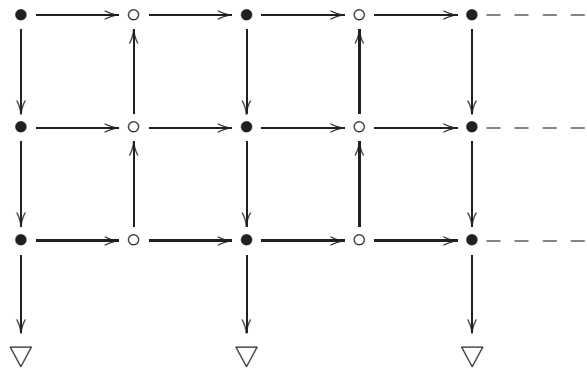
\end{document}